\documentclass{article}

\newcommand{\sfg}{{\sf g}}

\newcommand{\sfj}{{\sf j}}

\newcommand{\sfZF}{{\sf ZF}}
\newcommand{\sfZFL}{{\sf ZFL}}

\newcommand{\sfqk}{{\sf qk}}

\usepackage{proof}
\usepackage{latexsym}
\usepackage{amsfonts}
\usepackage{amssymb}
\usepackage{cite}

\newcommand{\blem}{\begin{lemma}}
\newcommand{\elem}{\end{lemma}}
\newcommand{\bth}{\begin{theorem}}
\newcommand{\ethm}{\end{theorem}}
\newcommand{\benu}{\begin{enumerate}}
\newcommand{\eenu}{\end{enumerate}}
\newcommand{\bdes}{\begin{description}}
\newcommand{\edes}{\end{description}}
\newcommand{\bdf}{\begin{definition}}
\newcommand{\edf}{\end{definition}}
\newcommand{\bcor}{\begin{cor}}
\newcommand{\ecor}{\end{cor}}
\newcommand{\bprp}{\begin{proposition}}
\newcommand{\eprp}{\end{proposition}}
\newcommand{\bmlem}{\begin{mlemma}}
\newcommand{\emlem}{\end{mlemma}}
\newcommand{\bclm}{\begin{claim}}
\newcommand{\eclm}{\end{claim}}
\newcommand{\bprf}{{\bf Proof}.\hspace{2mm}}

\newcommand{\eprf}{\hspace*{\fill} $\Box$}

\newcommand{\beqn}{\begin{equation}}
\newcommand{\eeqn}{\end{equation}}
\newcommand{\beqnarr}{\begin{eqnarray}}
\newcommand{\eeqnarr}{\end{eqnarray}}
\newcommand{\beqnarrs}{\begin{eqnarray*}}
\newcommand{\eeqnarrs}{\end{eqnarray*}}

\newcommand{\smlskp}{\\ \smallskip \\}
\newcommand{\spand}{\,\&\,}

\newcommand{\bfG}{\mbox{\boldmath$G$} }

\newcommand{\restrict}{\!\upharpoonright\!}

\newtheorem{theorem}{Theorem}[section]
\newtheorem{definition}[theorem]{Definition}
\newtheorem{proposition}[theorem]{Proposition}
\newtheorem{lemma}[theorem]{Lemma}
\newtheorem{cor}[theorem]{Corollary}
\newtheorem{mlemma}[theorem]{Main Lemma}
\newtheorem{claim}[theorem]{Claim}
\newtheorem{cyantheorem}{Theorem}

\newcommand{\alp}{\alpha}
\newcommand{\eps}{\epsilon}
\newcommand{\veps}{\varepsilon}
\newcommand{\del}{\delta}
\newcommand{\Del}{\Delta}
\newcommand{\ome}{\omega}
\newcommand{\Ome}{\Omega}
\newcommand{\bet}{\beta}
\newcommand{\gam}{\gamma}
\newcommand{\Gam}{\Gamma}
\newcommand{\kap}{\kappa}
\newcommand{\sig}{\sigma}
\newcommand{\Sig}{\Sigma}
\newcommand{\tht}{\theta}
\newcommand{\Tht}{\Theta}
\newcommand{\lam}{\lambda}
\newcommand{\Lam}{\Lambda}
\newcommand{\vphi}{\varphi}

\newcommand{\fal}{\forall}
\newcommand{\exi}{\exists}

\newcommand{\Rarw }{\Rightarrow}

\newcommand{\lrarw}{\leftrightarrow}
\newcommand{\Lrarw}{\Leftrightarrow}

\newcommand{\darw}{\downarrow}

\newcommand{\cala}{{\cal A}}

\newcommand{\cald}{{\cal D}}

\newcommand{\calf}{{\cal F}}
\newcommand{\calg}{{\cal G}}
\newcommand{\calh}{{\cal H}}
\newcommand{\cali}{{\cal I}}

\newcommand{\calt}{{\cal T}}

\newcommand{\calL}{{\cal L}}
\newcommand{\calk}{{\cal K}}

\newcommand{\calP}{{\cal P}}

\newcommand{\calx}{{\cal X}}

\newcommand{\st}{\star}

\newcommand{\sfv}{{\sf v}}
\newcommand{\sfe}{{\sf e}}

\newcommand{\sfk}{{\sf k}}
\newcommand{\rk}{\mbox{{\rm rk}}}

\newcommand{\la}{\langle}
\newcommand{\ra}{\rangle}

\newcommand{\setm}{\setminus}

\title{
Lifting up the proof theory to the countables
: Zermelo-Fraenkel set theory
}
\author{Toshiyasu Arai
\\
Graduate School of Science,
Chiba University
\\
1-33, Yayoi-cho, Inage-ku,
Chiba, 263-8522, JAPAN
\\
tosarai@faculty.chiba-u.jp
}
\date{}
\begin{document}
\maketitle

\begin{abstract}
We describe the countable ordinals in terms of iterations of Mostowski collapsings. 
This gives a proof-theoretic bound on definable countable ordinals
in Zermelo-Fraenkel set theory $\sfZF$.
\end{abstract}

\section{Introduction}

In these decades ordinal analyses (mainly of set theories) have progressed greatly, 
cf. M. Rathjen's contributions \cite{Rathjen91, Rathjen94, RathjenAFML1, RathjenAFML2} and 
\cite{ptMahlo, ptpi3, ptgentzen}.

Current ordinal analyses are {\it recursive\/}.
By recursive ordinal analyses we mean that everything in the analyses is recursive (on $\ome$).
Namely notation systems for ordinals to measure the proof-theoretic strengths of formal theories
are recursive, and operations on (codes of recursive) infinite derivations to eliminate cut inferences
are recursive, and so on.
Moreover in the analyses we consider only derivations of recursive statements on the least recursively regular ordinal $\ome_{1}^{CK}$.
We now ask: Can we lift up recursive ordinal analyses to countables
through a non-effective ordinal analysis?
By an analysis on countables we aim at bounding provability in formal theories for sets with respect to
statements on countable sets.

The proof technique in these ordinal analyses (cut-elimination with collapsing functions) has been successful in
describing the bounds on provability in
theories on {\it recursive analogues\/}
of (small) large cardinals, which were introduced by Richter and Aczel\cite{Richter-Aczel74}.
We can expect that the technique works also for set theories of (true) large cardinals.
In this paper we give a way to describe a bound on provability
in Zermelo-Fraenkel set theory $\sfZF$.
We describe the countable ordinal $\Psi_{\ome_{1}}\veps_{I+1}$, and show that
the ordinal is a proof-theoretic bound on definable countable ordinals 
provably existing
in Zermelo-Fraenkel set theory $\sfZF$, Theorem \ref{th:mainthZ}.

Let us describe the content of this paper. 
In section \ref{sect:Skolemhull} we give a characterization of the regularity of ordinals
in terms of $\Sig_{1}$-Skolem hulls.
In section \ref{sect:Ztheory} we introduce a theory of sets which
is equivalent to $\sfZF+(V=L)$, 
and in section \ref{sect:ordinalnotation} 
collapsing functions $\alp\mapsto\Psi_{\kap,n}\alp<\kap$ are introduced 
for each uncountable regular cardinal $\kap\leq I$ and $n<\ome$, cf. Definition \ref{df:Cpsiregularsm},
where $I$ is intended to denote the least weakly inaccessible cardinal.
Let $\ome_{k}(I+1)$ denote the tower of $\ome$ with the next epsilon number 
$\veps_{I+1}=\sup\{\ome_{k}(I+1): k<\ome\}$ above $I$.
It is easy to see that the predicate $x=\Psi_{\kap,n}\alp$ is a $\Sig_{n+1}$-predicate for $\alp<\veps_{I+1}$, and
for each $n,k<\ome$
$\sfZF+(V=L)$ proves $\fal\alp<\ome_{k}(I+1)
\fal\kap\leq I\exi x<\kap[x=\Psi_{\kap,n}\alp]
$, cf. Lemma \ref{lem:lowerbndreg}.

Conversely 
we show the following Theorem \ref{th:mainthZ} in the fragment $I\Sig^{0}_{1}$ of first-order arithmetic.

\bth\label{th:mainthZ}

For a sentence $\exi x\in L_{\ome_{1}}\, \vphi(x)$ with a first-order formula $\vphi(x)$, 
if
\[
\sfZF+(V=L)\vdash\exi x\in L_{\ome_{1}}\,\vphi(x) 
\]
then
\[
\exi n<\ome[\sfZF+(V=L)\vdash \exi x\in L_{\Psi_{\ome_{1},n}\ome_{n}(I+1)} \vphi(x)]
.\]
\end{theorem}
{\bf Remark}.
From Theorem \ref{th:mainthZ} together with  Lemma \ref{lem:lowerbndreg}
it follows that the countable ordinal
\[
\Psi_{\ome_{1}}\veps_{I+1}:=\sup\{\Psi_{\ome_{1},n}\ome_{n}(I+1): n<\ome\}
\]
is the limit of $\sfZF+(V=L)$-provably countable ordinals in the following sense:
\[
\Psi_{\ome_{1}}\veps_{I+1} = \sup\{\alp<\ome_{1}: \alp \mbox{ is a }
\sfZF+(V=L) \mbox{-provably countable ordinal }\}
\]
where by saying that an ordinal $\alp$ is a $\sfZF+(V=L)$-provably countable we mean
\[
\sfZF+(V=L)\vdash \exi! x<\ome_{1}\,\vphi(x) \spand L\models\vphi(\alp) \mbox{ for some formula } \vphi
.\]

From Theorem \ref{th:mainthZ} we see that if $\sfZF+(V=L)$ proves the existence of a real $a\in{}^{\ome}\ome$ 
enjoying a first-order condition $\vphi(a)$, 
$\sfZF+(V=L)\vdash\exi a\in{}^{\ome}\ome \vphi(a)$,
then such a real $a$ is already in level $L_{\Psi_{\ome_{1}}\veps_{I+1}}$ of constructible hierarchy.
\\

This paper is based on a technique, {\it operator controlled derivations\/},
which was introduced by W. Buchholz\cite{Buchholz},
hereby he gave a convincing ordinal analysis for the theory $\mbox{{\rm KP}}i$ of recursively inaccessible ordinals, which is a recursive analogue of $\sfZF$.
In section \ref{sect:controlledOme} operator controlled derivations for $\sfZF$
are introduced, and in the final section \ref{sect:proof}
Theorem \ref{th:mainthZ}  is concluded.
First let us explain the technique briefly.

In an operator controlled derivation, ordinals occurring in the derivation are controlled by an operator
$\calh$ on ordinals.
Through this we see that these ordinals are in a Skolem hull $\calh$.
On the other side a recursive notation system is defined through an iteration of Skolem hullings.
Suppose that a formal theory on sets proves a sentence $\exi x<\ome_{1}^{CK}\tht$ for a bounded formula $\tht$.
Then the technique tells us how many times do we iterate Skolem hulllings 
to bound a recursive ordinal $x$, a witness for $\tht$.

To be specific, let us explain how a Skolem hull looks like.
Let $\calf$ denote a set of functions.

\bdf\label{df:clsgeneral}{\rm (Cf. \cite{BuchholzFeferman}.)}
{\rm For sets} $X$,
$Cl(X;\calf)$ {\rm denotes the} Skolem hull {\rm of} $X$ {\rm under the functions in} $\calf$.

{\rm The set} $Cl(X;\calf)$ {\rm is inductively generated as follows.}
\benu
\item 
$X\subset Cl(X;\calf)$.
\item
{\rm If} $\vec{x}\subset Cl(X;\calf)$, $f\in\calf$ {\rm and} $\vec{x}\subset dom(f)${\rm , then} $f(\vec{x})\in Cl(X;\calf)$.
\eenu
\edf

Now let us restrict the construction on the class of ordinals $Ord$.
$\Ome=\ome_{1}$ denotes the least uncountable ordinal.
Let $\calf$ be a {\it countable\/} set of ordinal functions $f:Ord^{n}\to Ord$, where
the arity $n<\ome$ of the function $f$ is fixed for each $f$.
Assume that 0-ary functions $0,\Ome$ belong to $\calf$.

\bprp\label{prp:closedpnt}
\benu
\item\label{prp:closedpnt.0}
$
\fal\alp<\Ome\exi\bet<\Ome[Cl(\alp;\calf)\cap\Ome\subset\bet]
$.

\item\label{prp:closedpnt.1}
$
\fal\alp<\Ome\exi\bet<\Ome[\bet>\alp\spand Cl(\bet;\calf)\cap\Ome\subset\bet]
$.
Namely $\{\bet<\Ome: Cl(\bet;\calf)\cap\Ome\subset\bet\}$ is unbounded in $\Ome$.

\item\label{prp:closedpnt.2}
$\{\bet<\Ome: Cl(\bet;\calf)\cap\Ome\subset\bet\}$ is closed in $\Ome$.
\eenu
\eprp
\bprf
\ref{prp:closedpnt}.\ref{prp:closedpnt.0}.
If $\alp<\Ome$, then the set $Cl(\alp;\calf)$ is countable.
\\
\ref{prp:closedpnt}.\ref{prp:closedpnt.1}.
Given $\alp<\Ome$, define $\{\bet_{n}\}_{n}$ inductively,
$\bet_{0}=\alp+1$, $\bet_{n+1}=\min\{\bet<\Ome: Cl(\bet_{n};\calf)\cap\Ome\subset\bet\}$.
Then $\bet=\sup_{n}\bet_{n}$ is a desired one.
$\bet<\Ome$ since $\Ome$ is regular.
\eprf
\\
Let us enumerate the closed points.
Define sets $Cl_{\alp}(X;\calf)$ and ordinals $\psi_{\Ome}(\alp;\calf)$ by simultaneous recursion on ordinals $\alp$
as follows.

Let 
\[
Cl_{\alp}(X;\calf):=Cl(X;\calf\cup\{\psi_{\Ome}(\cdot;\calf)\restrict\alp\})
\]
where
\[
\psi_{\Ome}(\alp;\calf)=\min\{\bet\leq\Ome: Cl_{\alp}(\bet;\calf)\cap\Ome\subset\bet\}
.\]
Then a transfinite induction on $\alp$ shows with Proposition \ref{prp:closedpnt}.\ref{prp:closedpnt.1}
\[
\fal\alp\exi\bet<\Ome[\psi_{\Ome}(\alp;\calf)=\bet]
.\]

For $\calf_{0}=\{0,\Ome\}\cup\{\lam xy. x+y, \lam x. \ome^{x}\}$(and the Veblen function $\lam xy.\vphi xy$),
$\psi_{\Ome}(\veps_{\Ome+1} ;\calf_{0})$ is the Howard ordinal, the proof-theoretic ordinal of
the theory $\mbox{ID}_{1}$ for non-iterated positive elementary inductive definition on $\ome$,
or equivalently of $\mbox{KP}\ome$, i.e., 
Kripke-Platek set theory with the axiom of infinity.

Observe that each function in $\calf_{0}$ is $\{\Ome\}$-recursive in $L_{\sig}$ for any $\sig>\Ome$.
Here an $\{\Ome\}$-recursive function is $\Sig$-definable from the 
0-ary function $\Ome$, a parameter.

Now let us extend $\calf_{0}$ to the set $\calf_{all}$
 of {\it all\/} $\{\Ome\}$-recursive functions on $L_{\sig}$.
Then it turns out that 
 $Cl(\alp;\calf_{all})$ is the $\Sig_{1}$-Skolem hull $\mbox{{\rm Hull}}_{\Sig_{1}}^{\sig}(\alp\cup\{\Ome\})$
of $\alp\cup\{\Ome\}$ on $L_{\sig}$, and this gives a characterization of the regularity of the ordinal $\Ome$,
cf. Theorem \ref{th:cofinalitylocal} below.


\section{$\Sig_{n}$-Skolem hulls}\label{sect:Skolemhull}

For a model $\la M;\in\restrict (M\times M)\ra$ and $X\subset M$,
$\Sig_{n}^{M}(X)$ denotes the set of $\Sig_{n}(X)$-definable subsets of $M$,
where $\Sig_{n}(X)$-formulae may have parameters from $X$.
$\Sig_{n}^{M}(M)$ is denoted $\Sig_{n}(M)$.

An ordinal $\alp>1$ is said to be a {\it multiplicative principal number\/} iff
$\alp$ is closed under ordinal multiplication, i.e., $\exi \bet[\alp=\ome^{\ome^{\bet}}]$.
If $\alp$ is a multiplicative principal number, then $\alp$ is closed under G\"odel's pairing function $j$
and there exists a $\Del_{1}$-bijection between $\alp$ and $L_{\alp}$
for the constructible hierarchy $L_{\alp}$ up to $\alp$.
In this section 
$\sig$ is assumed to be a multiplicative principal number$>\ome$.

\bdf\label{df:crdS} 

\benu
\item
$Reg$ {\rm denotes the class of uncountable regular ordinals.}

\item 
$
cf(\kap) := \min\{\alp\leq\kap : \mbox{{\rm there is a cofinal map }} f:\alp\to \kap\}
$.
\beqnarrs
\kap \mbox{ {\rm is} } uncountable \: regular  & :\Lrarw & \kap\in Reg
\Lrarw  \kap>\ome\spand cf(\kap)=\kap 
\\
&  \Lrarw & \kap>\ome\spand \fal \alp<\kap(\alp<cf(\kap))
\eeqnarrs
$
 card(\alp)<card(\kap) :\Lrarw  \mbox{{\rm there is no surjective map }} f:\alp\to\kap
$.


\item
$\rho(L_{\sig})$ {\rm denotes the} $\Sig_{1}$-projectum {\rm of} $L_{\sig}${\rm :}
$\rho(L_{\sig})$ {\rm is the least ordinal} $\rho$ {\rm such that}
$
\calP(\rho)\cap\Sig_{1}(L_{\sig})\not\subset L_{\sig}
$.

\item
{\rm Let} $\alp\leq\bet$ {\rm and} $f: L_{\alp}\to L_{\bet}$.
{\rm Then the map} $f$ {\rm is a} $\Sig_{n}$-elementary embedding{\rm , denoted}
$
f: L_{\alp}\prec_{\Sig_{n}} L_{\bet}
$
{\rm iff for any} $\Sig_{n}(L_{\alp})${\rm -sentence} $\vphi[\bar{a}]\, (\bar{a}\subset L_{\alp})$,
$
L_{\alp} \models \vphi[\bar{a}] \Lrarw 
L_{\bet} \models \vphi[f(\bar{a})]
$
{\rm where} $f(\bar{a})=f(a_{1}),\ldots,f(a_{k})$ {\rm for} $\bar{a}=a_{1},\ldots,a_{k}$.
{\rm An ordinal} $\gam$ {\rm such that}
$
\fal\del<\gam[f(\del)=\del]\spand f(\gam)>\gam
$
{\rm is said to be the} critical point {\rm of the} $\Sig_{n}$-
{\rm elementary embedding} $f$
{\rm if such an ordinal} $\gam$ {\rm exists.}

\item
{\rm For} $X\subset L_{\sig}$,
$\mbox{{\rm Hull}}^{\sig}_{\Sig_{n}}(X)$ {\rm denotes the set} 
{\rm (}$\Sig_{n}$-Skolem hull {\rm of} $X$ {\rm in} $L_{\sig}${\rm ) defined as follows.}
$<_{L}$ {\rm denotes a} $\Del_{1}${\rm -well ordering of the constructible universe} $L$.
{\rm Let} $\{\vphi_{i}:i\in\ome\}$ {\rm denote an enumeration of} $\Sig_{n}${\rm -formulae in the language}
$\{\in\}${\rm . Each is of the form} $\vphi_{i}\equiv\exi y\theta_{i}(x,y;u)\, (\tht\in\Pi_{n-1})$ {\rm with fixed variables} $x,y,u${\rm . Set for} $b\in X$
\beqnarr
r_{\Sig_{n}}^{\sig}(i,b) & \simeq & \mbox{ {\rm the }} <_{L} \mbox{{\rm -least }} c\in L_{\sig}
\mbox{ {\rm such that} } L_{\sig}\models\theta_{i}((c)_{0},(c)_{1}; b)
\nonumber
\\
h_{\Sig_{n}}^{\sig}(i,b) & \simeq & (r_{\Sig_{n}}^{\sig}(i,b))_{0}
\label{eq:sighull}
\\
\mbox{{\rm Hull}}^{\sig}_{\Sig_{n}}(X) & = & rng(h_{\Sig_{n}}^{\sig})=\{h_{\Sig_{n}}^{\sig}(i,b)\in L_{\sig}:i\in\ome, b\in X\}
\nonumber
\eeqnarr
{\rm Then}
$
L_{\sig}\models \exi x\exi y\, \theta_{i}(x,y;b) \to h_{\Sig_{n}}^{\sig}(i,b)\darw \spand \exi y\, \theta_{i}(h_{\Sig_{n}}^{\sig}(i,b),y;b)
$.
\eenu
\edf

The following Propositions \ref{prp:hullcl}, \ref{prp:rhHull} and \ref{clm:crdlocal3.1} are easy to see.

\bprp\label{prp:hullcl}
For $a,\kap\in L_{\sig}$,
$
\mbox{{\rm Hull}}^{\sig}_{\Sig_{1}}(a\cup\{\kap\})=Cl(a;\calf_{all})
$,
where in the RHS, Definition \ref{df:clsgeneral}, $\Ome$ is replaced by 
$\kap\in\calf_{all}$, and
$\calf_{all}$ denotes the set of all $\{\kap\}$-recursive (partial) functions on $L_{\sig}$.
Namely $f\in\calf_{all}$ iff there exists an $i<\ome$ such that
$f(b)\simeq\bet\Lrarw  h_{\Sig_{1}}^{\sig}(i,\la b, \kap\ra)\simeq\bet$ for $b<\sig$,
where $\la b,c\ra$ denotes the pairing of $b$ and $c$.
\eprp



\bprp\label{prp:rhHull}
Assume that $X$ is a set in $L_{\sig}$.
Then
$r_{\Sig_{n}}^{\sig}$ and $h_{\Sig_{n}}^{\sig}$ are partial $\Del_{n}(L_{\sig})$-maps such that the domain of $h_{\Sig_{n}}^{\sig}$ is a 
$\Sig_{n}(L_{\sig})$-subset of 
$\ome\times X$.
Therefore its range $\mbox{{\rm Hull}}^{\sig}_{\Sig_{n}}(X)$ is a $\Sig_{n}(L_{\sig})$-subset of $L_{\sig}$.
%
\eprp


\bprp\label{clm:crdlocal3.1} 
Let $Y=\mbox{{\rm Hull}}^{\sig}_{\Sig_{n}}(X)$.
For any $\Sig_{n}(Y)$-sentence $\vphi(\bar{a})$ with parameters $\bar{a}$ from $Y$
$
L_{\sig}\models \vphi(\bar{a})\Lrarw Y\models \vphi(\bar{a})
$.
Namely
$
Y \prec_{\Sig_{n}}L_{\sig}
$.
\eprp

\bdf\label{df:pikap}{\rm (Mostowski collapsing function} $F${\rm )}

{\rm Let} $n\geq 1$.
{\rm By Proposition \ref{clm:crdlocal3.1} and the Condensation Lemma, cf.\,\cite{Devlin},
we have an isomorphism (Mostowski collapsing function) }
\[
F:\mbox{{\rm Hull}}^{\sig}_{\Sig_{n}}(X)\lrarw  L_{\gam}
\]
{\rm for an ordinal} $\gam\leq\sig$ {\rm such that} $F\restrict Y=id\restrict Y$ 
{\rm for any transitive}
$Y\subset \mbox{{\rm Hull}}^{\sig}_{\Sig_{n}}(X)$.

{\rm Let us denote, though} $\sig\not\in dom(F)=\mbox{{\rm Hull}}^{\sig}_{\Sig_{n}}(X)$
\[
F(\sig):=\gam
.\]
{\rm Also for the above Mostowski collapsing map} $F$ {\rm let}
\[
F^{\Sig_{n}}(x;\sig,X):=F(x)
.\]
{\rm The inverse} $G:=F^{-1}$ {\rm of} $F$ {\rm is a} $\Sig_{n}${\rm -elementary embedding
from} $L_{F(\sig)}$ {\rm to} $L_{\sig}$.
\edf

\bdf
{\rm Let} $\kap$ {\rm be an ordinal such that} $\ome<\kap<\sig${\rm , and let}
\[
F_{\bet\cup\{\kap\}}^{\Sig_{n}}(x):=F^{\Sig_{n}}(x;\sig,\bet\cup\{\kap\})
.\]
{\rm Then put}
\beqnarrs
C^{\sig}_{\Sig_{1}}(\kap) & := &\{x<\kap: x\in Cr^{\sig}_{\Sig_{1}}(\{\kap\}) \spand F_{x\cup\{\kap\}}^{\Sig_{1}}(\sig)<\kap\}
\\
x\in Cr^{\sig}_{\Sig_{1}}(\{\kap\}) & :\Lrarw & \mbox{{\rm Hull}}^{\sig}_{\Sig_{1}}(x\cup\{\kap\})\cap \kap\subset x
\eeqnarrs
\edf


\bprp\label{prp:sigprojectum}
Let $\alp$ be a multiplicative principal number with $\ome\leq\alp<\kap<\sig$.
Assume that $\sig$ is recursively regular 
and the $\Sig_{1}$-projectum $\rho(L_{\sig})>\alp$.

Then for the map $h_{\Sig_{1}}^{\sig}$ with $X=\alp\cup\{\kap\}$ in (\ref{eq:sighull}) we have $dom(h_{\Sig_{1}}^{\sig})\in L_{\sig}$.
Therefore $\mbox{{\rm Hull}}^{\sig}_{\Sig_{1}}(\alp\cup\{\kap\})=rng(h_{\Sig_{1}}^{\sig})$ is a set in $L_{\sig}$, 
and the Mostowski collapsing function 
$F_{\alp\cup\{\kap\}}^{\Sig_{1}}:\mbox{{\rm Hull}}^{\sig}_{\Sig_{1}}(\alp\cup\{\kap\}) \lrarw L_{F^{\Sig_{1}}_{\alp\cup\{\kap\}}(\sig)}$
is a $\Del_{1}(L_{\sig})$-map.
Hence $L_{F^{\Sig_{1}}_{\alp\cup\{\kap\}}(\sig)}=rng(F^{\Sig_{1}}_{\alp\cup\{\kap\}})\in L_{\sig}$, i.e., $F^{\Sig_{1}}_{\alp\cup\{\kap\}}(\sig)<\sig$.

Moreover if $\rho(L_{\sig})>\kap$, then
$Cr^{\sig}_{\Sig_{1}}(\{\kap\})=\{x<\kap: \mbox{{\rm Hull}}^{\sig}_{\Sig_{1}}(x\cup\{\kap\})\cap \kap\subset x\}$ 
is a set in $L_{\sig}$.
\eprp
\bprf
By the definition
$dom(h_{\Sig_{1}}^{\sig})=\{(i,\bet)\in\ome\times\alp: L_{\sig}\models\exi c \, \theta_{i}((c)_{0},(c)_{1};\bet,\kap)\}$
is a $\Sig_{1}(L_{\sig})$-subset of $\ome\times\alp\lrarw\alp$. 

By the supposition
we have $\alp<\rho(L_{\sig})$.
Therefore any $\Sig_{1}(L_{\sig})$-subset of $\alp$ is a set in $L_{\sig}$ by the definition of the $\Sig_{1}$-projectum.

$Cr^{\sig}_{\Sig_{1}}(\{\kap\})$ is a $\Pi_{1}(L_{\sig})$-subset of $\kap<\rho(L_{\sig})$, and hence is a set in $L_{\sig}$.
\eprf

\blem\label{lem:cofinalitylocal}
Let $\alp$ be a multiplicative principal number with $\ome\leq\alp<\kap<\sig$.
Assume that $\sig$ is recursively regular and $L_{\sig}\models \alp<cf(\kap)$.
\benu
\item\label{lem:cofinalitylocal1}
$\alp<\rho(L_{\sig})$.

\item\label{lem:cofinalitylocal15}
$F^{\Sig_{1}}_{\alp\cup\{\kap\}}(\sig)<\kap$.

\item\label{lem:cofinalitylocal2}
Let $\bet$ denote the least ordinal $\bet\leq\kap$ such that
$\mbox{{\rm Hull}}^{\sig}_{\Sig_{1}}(\alp\cup\{\kap\})\cap\kap\subset\bet$.
Then $\bet<\kap$ and  $L_{\sig}\models \bet<cf(\kap)$, and hence $\bet<\rho(L_{\sig})$.
\eenu
\elem
\bprf
\\
\ref{lem:cofinalitylocal}.\ref{lem:cofinalitylocal1}(Cf. \cite{Barwise}.).
 Let $\emptyset\neq B\in\Sig_{1}(L_{\sig})\cap\calP(\alp)$.
 We show $B\in L_{\sig}$.
Let $g:\sig\to B$ be a surjection, 
and $f$ be the map
$
f(\gam)=g(\mu \del(g(\del)\not\in \{f(\xi):\xi<\gam\}))
$,
 i.e., $f(\gam)$ is the $\gam$th member of $B$.
Both $g$ and $f$ are $\Del_{1}(L_{\sig})$-maps.
Suppose that $f$ is total.
The $\Sig_{1}(L_{\sig})$-injection $f$ from $\sig$ to $\alp$ yields an injection from $\kap$ to $\alp$ in $L_{\sig}$,
whose inverse would be a cofinal map from $\alp$ to $\kap$ in $L_{\sig}$.
Let $\gam_{0}$ be the least $\gam<\sig$ such that $f(\gam)$ is undefined.
Then $B=\{f(\gam):\gam<\gam_{0}\}$, and hence $B\in L_{\sig}$ by $\Sig$-Replacement.
\\
\ref{lem:cofinalitylocal}.\ref{lem:cofinalitylocal15}.
We have
$
\alp<\rho(L_{\sig})
$ by Lemma \ref{lem:cofinalitylocal}.\ref{lem:cofinalitylocal1}.
Then by Proposition \ref{prp:sigprojectum}  we have
$F^{\Sig_{1}}_{\alp\cup\{\kap\}}(\sig)<\kap$.
\\
\ref{lem:cofinalitylocal}.\ref{lem:cofinalitylocal2}.
By Proposition \ref{prp:sigprojectum}, there exists a surjective map in $L_{\sig}$ from $\alp$
to $\mbox{{\rm Hull}}^{\sig}_{\Sig_{1}}(\alp\cup\{\kap\})$.
Therefore
$\mbox{{\rm Hull}}^{\sig}_{\Sig_{1}}(\alp\cup\{\kap\})\cap\kap$ is bounded in $\kap$.
By the minimality of $\bet$, $\mbox{{\rm Hull}}^{\sig}_{\Sig_{1}}(\alp\cup\{\kap\})\cap\kap$ is cofinal in $\bet$.
\eprf

\bprp\label{prp:complexcr}
Let $n\geq 1$ and
$L_{\sig}\models\mbox{{\rm KP}}\ome+\Sig_{n}\mbox{{\rm -Collection}}$.
Then for $\kap\leq\sig$,
$\{(x,y): x<\kap\spand y=\min\{y<\kap: \mbox{{\rm Hull}}^{\sig}_{\Sig_{n}}(x\cup\{\kap\})\cap\kap\subset y\}\}$ 
is a $Bool(\Sig_{n}(L_{\sig}))$-predicate on $\kap$, and hence a set in $L_{\sig}$
if $\kap<\sig$ and $L_{\sig}\models\Sig_{n}\mbox{{\rm -Separation}}$.
\eprp
\bprf
Let $\vphi(y,\kap)$ be the $\Pi_{n}$-predicate $\vphi(y,\kap):\Lrarw 
\fal z<\kap[z\in\mbox{{\rm Hull}}^{\sig}_{\Sig_{n}}(x\cup\{\kap\})\to z\in y]$.
Then
$y=\min\{y<\kap: \mbox{{\rm Hull}}^{\sig}_{\Sig_{n}}(x\cup\{\kap\})\cap\kap\subset y\}$ iff
$y<\kap\land \vphi(y,\kap)
\land \fal u<y\lnot \vphi(u,\kap)$, which is  
$Bool(\Sig_{n}(L_{\sig}))$ by $\Pi_{n-1}\mbox{{\rm -Collection}}$.
\eprf

\subsection{Regularity}

$F^{\Sig_{1}}_{x\cup\{\kap\}}(y)$ denotes the Mostowski collapse $F^{\Sig_{1}}(y;\sig,x\cup\{\kap\})$.
The following Theorems \ref{th:cofinalitylocal} and \ref{th:cofinalitylocalpower} should be folklore.

\bth\label{th:cofinalitylocal}{\rm (Cf. \cite{sneak}.)}
Let $\sig$ be an ordinal such that $L_{\sig}\models \mbox{{\rm KP}}\ome+\Sig_{1}\mbox{{\rm -Separation}}$,
 and $\ome\leq\alp<\kap<\sig$ with
$\alp$ a multiplicative principal number and $\kap$ a limit ordinal.
Then the following conditions are mutually equivalent:
\benu

\item 
\beqn\label{eqn:cofinalitylocal2}
L_{\sig}\models {}^{\alp}\kap\subset L_{\kap}
\eeqn

\item 
\beqn\label{eqn:cofinalitylocal3}
L_{\sig}\models \alp<cf(\kap)
\eeqn

\item
There exists an ordinal $x$ such that $x\in C^{\sig}_{\Sig_{1}}(\kap)\cap(\alp,\kap)$, i.e., 
\beqn\label{eqn:cofinalitylocal0}
x\in Cr^{\sig}_{\Sig_{1}}(\{\kap\})\cap(\alp,\kap)
\spand
F^{\Sig_{1}}_{x\cup\{\kap\}}(\sig)<\kap
\eeqn



\item 
For the Mostowski collapse
$F^{\Sig_{1}}_{x\cup\{\kap\}}(y)$
\beqnarr
&& \exi x[\alp<x=F^{\Sig_{1}}_{x\cup\{\kap\}}(\kap)<F^{\Sig_{1}}_{x\cup\{\kap\}}(\sig)<\kap\spand  \fal \Sig_{1}\,\vphi\fal a\in L_{x}
 \nonumber \\
&&(L_{\sig}\models\vphi[\kap,a]\to  L_{F^{\Sig_{1}}_{x\cup\{\kap\}}(\sig)} \models\vphi[x,a])] \label{eqn:cofinalitylocal1}
\eeqnarr



\eenu

\end{theorem}
\bprf
Obviously under the assumption that $\sig$ is recursively regular,
(\ref{eqn:cofinalitylocal2}) and (\ref{eqn:cofinalitylocal3}) are mutually equivalent, and 
(\ref{eqn:cofinalitylocal0}) implies (\ref{eqn:cofinalitylocal1}).

Assume $\sig$ is recursively regular, $\kap$ denotes a limit ordinal and $\alp$ a multiplicative principal number with $\ome\leq\alp<\kap<\sig$.
\\

\noindent
 (\ref{eqn:cofinalitylocal1}) $\Rarw $ (\ref{eqn:cofinalitylocal2}).
Suppose there exist an ordinal $x$
such that
$\alp<x= F^{\Sig_{1}}_{x\cup\{\kap\}}(\kap)<F^{\Sig_{1}}_{x\cup\{\kap\}}(\sig)<\kap$ and for any $\Sig_{1}$ $\vphi$ and
any $a\in L_{x}$
\beqn
\renewcommand{\theequation}{\ref{eqn:cofinalitylocal1}}
L_{\sig}\models\vphi[\kap,a]\Rarw  L_{F^{\Sig_{1}}_{x\cup\{\kap\}}(\sig)}\models\vphi[x,a]
\eeqn
\addtocounter{equation}{-1}
Let us show 
\[
L_{\sig}\models {}^{\alp}\kap\subset L_{\kap}
.\]
Define a $\Del_{1}(L_{\sig})$-partial map $S:dom(S)\to{}^{\alp}\kap\cap L_{\kap}\, (dom(S)\subset\kap)$ by
letting $S_{\bet}$ be the $<_{L}$ least $X\in{}^{\alp}\kap\cap L_{\kap}$ such that $\fal\gam<\bet(X\neq S_{\gam})$. 

It suffices to show that $L_{\sig}\models {}^{\alp}\kap\subset\{S_{\bet}\}_{\bet}=rng(S)$. 
Suppose there exists an $f\in{}^{\alp}\kap\cap L_{\sig}$ so that $\fal\bet<\kap(S_{\bet}\neq f)$ and let $f_{0}$ denote the $<_{L}$-least such function.
Then $f_{0}$ is $\Sig_{1}$ definable on $L_{\sig}$ from $\{\alp,\kap\}$: for the $\Del_{1}(L_{\sig})$-formula
$
\vphi(f,\alp,\kap) :\Lrarw \theta(f,\alp,\kap)\land \fal g<_{L}f \, \lnot\theta(g,\alp,\kap)
$
with
$
\theta(f,\alp,\kap) :\Lrarw f\in {}^{\alp}\kap \land\fal\bet<\kap(S_{\bet}\neq f)
$
we have
$
L_{\sig}\models\vphi(f_{0},\alp,\kap)\spand L_{\sig}\models\exi !f\vphi(f,\alp,\kap)
$.
By (\ref{eqn:cofinalitylocal1}) 
we have $L_{F^{\Sig_{1}}_{x\cup\{\kap\}}(\sig)}\models\exi f\vphi(f,\alp,x)$, i.e., there exists the $<_{L}$-least $f_{1}\in  {}^{\alp}x\cap L_{\kap}\,(F^{\Sig_{1}}_{x\cup\{\kap\}}(\sig)\leq\kap)$ such that $\fal\bet<x(<\kap)(S_{\bet}\neq f_{1})$. 

We show $L_{\kap}\ni f_{1}=f_{0}$. This yields a contradiction. 
It suffices to see $f_{1}\subset f_{0}$ for $f_{1}:\alp\to x$ and $f_{0}:\alp\to\kap$.
By (\ref{eqn:cofinalitylocal1}) 
we have for $\bet<\alp$, $\gam<x$
\beqnarrs
&& f_{1}(\bet)=\gam \Lrarw  L_{F^{\Sig_{1}}_{x\cup\{\kap\}}(\sig)}\models\fal f[\vphi(f,\alp,x)\to f(\bet)=\gam] \Rarw \\
&& L_{\sig}\models\fal f[\vphi(f,\alp,\kap)\to f(\bet)=\gam] \Lrarw f_{0}(\bet)=\gam
\eeqnarrs

Note that in this proof it suffices to assume that $\sig$ is recursively regular,
and we see that the condition $F^{\Sig_{1}}_{x\cup\{\kap\}}(\sig)<\kap$ can be weakened to
$F^{\Sig_{1}}_{x\cup\{\kap\}}(\sig)\leq\kap$ in (\ref{eqn:cofinalitylocal0}) and (\ref{eqn:cofinalitylocal1}).
\\
(\ref{eqn:cofinalitylocal3}) $\Rarw$ (\ref{eqn:cofinalitylocal0}). 
Assume $L_{\sig}\models\Sig_{1}\mbox{{\rm -Separation}}$, and $L_{\sig}\models \alp<cf(\kap)$.






We show the existence of an ordinal $x<\kap$ such that
\[
x>\alp \spand \mbox{{\rm Hull}}^{\sig}_{\Sig_{1}}(x\cup\{\kap\})\cap\kap\subset x
\spand F^{\Sig_{1}}_{x\cup\{\kap\}}(\sig)<\kap
.\]
Then $F^{\Sig_{1}}_{x\cup\{\kap\}}(\kap)=x$.

As in the proof of Proposition \ref{prp:closedpnt}.\ref{prp:closedpnt.1}, 
define recursively ordinals $\{x_{n}\}_{n}$ as follows.
$x_{0}=\alp+1$, and $x_{n+1}$ is defined to be the least ordinal $x_{n+1}\leq\kap$ such that
$\mbox{{\rm Hull}}^{\sig}_{\Sig_{1}}(x_{n}\cup\{\kap\})\cap\kap\subset x_{n+1}$.
We see inductively that $x_{n}<\kap$ from Lemma \ref{lem:cofinalitylocal}.\ref{lem:cofinalitylocal2}.
On the other hand we have 
${}^{n}\kap\subset L_{\kap}$ by (\ref{eqn:cofinalitylocal2}).
Moreover by Proposition \ref{prp:complexcr},
the map $n\mapsto x_{n}$ is 
a $\Del_{1}$-set in $L_{\sig}\models\Sig_{1}\mbox{{\rm -Separation}}$.


Therefore $x=\sup_{n}x_{n}<\kap$ enjoys 
$x>\alp$, and 
$
\mbox{{\rm Hull}}^{\sig}_{\Sig_{1}}(x\cup\{\kap\})\cap\kap\subset x
$.

It remains to see $F^{\Sig_{1}}_{x\cup\{\kap\}}(\sig)<\kap$.
By Lemma \ref{lem:cofinalitylocal}.\ref{lem:cofinalitylocal15} it suffices to see
$x<cf(\kap)$.

Since there exists a $\Del_{1}(L_{\sig})$-surjective map 
$h_{n}:x_{n}\to \mbox{{\rm Hull}}^{\sig}_{\Sig_{1}}(x_{n}\cup\{\kap\})$,
pick an increasing cofinal map $f_{n}:x_{n}\to x_{n+1}$ in $L_{\sig}$
using the minimality of $x_{n+1}$.
Using the uniformity of $f_{n}$, we see the existence of an increasing cofinal map
$f:\alp\to x$ in $L_{\sig}$.
Therefore $L_{\sig}\models x<cf(\kap)$.
%
\eprf
\\

\noindent
{\bf Remark}.
In the proof of Theorem \ref{th:cofinalitylocal}, the assumption that
$L_{\sig}\models\Sig_{1}\mbox{{\rm -Separation}}$
 is used only in the part
(\ref{eqn:cofinalitylocal3}) $\Rarw$ (\ref{eqn:cofinalitylocal0}), and
everything except the part 
holds when $\sig$ is recursively regular.


\bcor\label{cor:cofinalreg}
Suppose $\kap$ is uncountable regular in $L_{\sig}\models\mbox{{\rm KP}}\ome+\Sig_{1}\mbox{{\rm -Separation}}$.
\benu
\item\label{cor:cofinalreg1}
$\kap$ is $\sig$-stable, i.e., $L_{\kap}\prec_{\Sig_{1}}L_{\sig}$.

\item\label{cor:cofinalreg15}
$\{\lam<\kap: \lam\in Reg\}=\{\lam<\kap: \lam \mbox{ {\rm is uncountable regular in }} L_{\sig}\}$
is a $\Del_{0}$-subset of $\kap$.
Therefore 
the map $\kap>\alp\mapsto\ome_{\alp}$ is a $\Del_{1}$-map on $L_{\sig}$.
On the other side the map $\sig>\alp\mapsto\ome_{\alp}$ is a $\Del_{2}$-map on $L_{\sig}$.


\eenu
\ecor
\bprf
\ref{cor:cofinalreg}.\ref{cor:cofinalreg1}.
Let $\vphi[a]$ be a $\Sig_{1}$-formula with a parameter $a\in L_{\kap}$.
Pick an $\alp_{a}\in C^{\sig}_{\Sig_{1}}(\kap)$ such that $a\in L_{\alp_{a}}$ by Theorem \ref{th:cofinalitylocal}.
Then 
$
L_{\sig}\models\vphi[a]\Rarw 
L_{F^{\Sig_{1}}_{\alp_{a}\cup\{\kap\}}}\models\vphi[a]
\Rarw
L_{\kap}\models\vphi[a]
$
 for $a=F^{\Sig_{1}}_{\alp_{a}\cup\{\kap\}}(a)$ and $F^{\Sig_{1}}_{\alp_{a}\cup\{\kap\}}(\sig)<\kap$.
\\
\ref{cor:cofinalreg}.\ref{cor:cofinalreg15}.
For $\lam<\kap$, we see from Corollary \ref{cor:cofinalreg}.\ref{cor:cofinalreg1},
$L_{\sig}\models \lam\in Reg \Lrarw L_{\kap}\models\lam\in Reg$.
\eprf
\\

For the existence of power sets we have the following Theorem \ref{th:cofinalitylocalpower}.

\bth\label{th:cofinalitylocalpower}{\rm (Cf. \cite{sneak}.)}
Let $\sig$ be recursively regular,
 and $\ome\leq\alp<\kap<\sig$ with
$\alp$ a multiplicative principal number and $\kap$ a limit ordinal.
Then the following conditions are mutually equivalent:
\benu

\item
\beqn\label{eqn:crdlocal-1}
\alp<\rho(L_{\sig})  \land 
F^{\Sig_{1}}_{\alp\cup\{\alp,\kap\}}(\sig)=F^{\Sig_{1}}(\sig;\sig,\alp\cup\{\alp,\kap\})<\kap
\eeqn

\item 
For the Mostowski collapse
$F^{\Sig_{1}}_{\alp\cup\{\alp,\kap\}}:\mbox{{\rm Hull}}_{\Sig_{1}}^{\sig}(\alp\cup\{\alp,\kap\})\lrarw L_{F^{\Sig_{1}}_{\alp\cup\{\alp,\kap\}}(\sig)}$
\beqnarr
&& \exi  x[\alp<x\leq F^{\Sig_{1}}_{\alp\cup\{\alp,\kap\}}(\kap)<F^{\Sig_{1}}_{\alp\cup\{\alp,\kap\}}(\sig)<\kap\spand  \fal \Sig_{1}\,\vphi\fal a\in L_{x}
 \nonumber \\
&&(L_{\sig}\models\vphi[\kap,a]\to  L_{F^{\Sig_{1}}_{\alp\cup\{\alp,\kap\}}(\sig)} \models\vphi[F^{\Sig_{1}}_{\alp\cup\{\alp,\kap\}}(\kap),a])] \label{eqn:crdlocal1}
\eeqnarr

\item 
\beqn\label{eqn:crdlocal2}
{\cal P}(\alp)\cap L_{\sig}\subset L_{\kap}
\eeqn

\item 
\beqn\label{eqn:crdlocal3}
L_{\sig}\models card(\alp)<card(\kap)
\eeqn

\eenu

\end{theorem}
\bprf
In showing the direction (\ref{eqn:crdlocal-1})$\Rarw$(\ref{eqn:crdlocal1}),
pick the least ordinal $x>\alp$ not in 
$\mbox{{\rm Hull}}^{\sig}_{\Sig_{1}}(\alp\cup\{\alp,\kap\})$.
(\ref{eqn:crdlocal2})$\Rarw$(\ref{eqn:crdlocal3}) 
and (\ref{eqn:crdlocal3})$\Rarw$(\ref{eqn:crdlocal-1}) 
are easily seen.
\\

\noindent
(\ref{eqn:crdlocal1})$\Rarw$(\ref{eqn:crdlocal2}).
As in the proof of  (\ref{eqn:cofinalitylocal1}) $\Rarw $ (\ref{eqn:cofinalitylocal2}),
define a $\Del_{1}$-partial map $S:dom(S)\to{\cal P}(\alp)\cap L_{\kap}\, (dom(S)\subset\kap)$ by letting
$S_{\bet}$ be the $<_{L}$ least $X\in{\cal P}(\alp)\cap L_{\kap}$ such that $\fal\gam<\bet(X\neq S_{\gam})$. 

It suffices to show that ${\cal P}(\alp)\cap L_{\sig}\subset\{S_{\bet}\}_{\bet}=rng(S)$. Suppose there exists an $X\in{\cal P}(\alp)\cap L_{\sig}$ so that $\fal\bet<\kap(S_{\bet}\neq X)$ and let $X_{0}$ denote the $<_{L}$-least such set. 
Then we see that $X_{0}$ is $\Sig_{1}$-definable in $L_{\sig}$ from $\{\alp,\kap\}$: there exists a $\Del_{1}$-formula $\vphi(X,\alp,\kap)$ such that
$
L_{\sig}\models\vphi(X_{0},\alp,\kap)\spand L_{\sig}\models\exi !X\vphi(X,\alp,\kap)
$.
By (\ref{eqn:crdlocal1}) we have 
$L_{F^{\Sig_{1}}_{\alp\cup\{\alp,\kap\}}(\sig)}\models\exi X\vphi(X,\alp,F^{\Sig_{1}}_{\alp\cup\{\alp,\kap\}}(\kap))$, 
i.e., there exists the $<_{L}$-least 
$X_{1}\in{\cal P}(\alp)\cap L_{F^{\Sig_{1}}_{\alp\cup\{\alp,\kap\}}(\sig)}\subset {\cal P}(\alp)\cap L_{\kap}$ 
such that $\fal\bet<F^{\Sig_{1}}_{\alp\cup\{\alp,\kap\}}(\kap)(<\kap)(S_{\bet}\neq X_{1})$. 
This means that $X_{1}=S_{F^{\Sig_{1}}_{\alp\cup\{\alp,\kap\}}(\kap)}$. 
We show $X_{1}=X_{0}$. This yields a contradiction. 
Denote $x\in a$ by $x\in^{+}a$ and $x\not\in a$ by $x\in^{-}a$. For any $\gam<\alp$, again by (\ref{eqn:crdlocal1}) we have
\beqnarrs
&& \gam\in^{\pm} X_{0}\Lrarw L_{\sig}\models\exi X(\gam\in^{\pm}X\land \vphi(X,\alp,\kap)) \Rarw \\
&& L_{F^{\Sig_{1}}_{\alp\cup\{\alp,\kap\}}(\sig)}\models\exi X(\gam\in^{\pm}X\land \vphi(X,\alp,F^{\Sig_{1}}_{\alp\cup\{\alp,\kap\}}(\kap))) \Lrarw \gam\in^{\pm}X_{1}
\eeqnarrs

\eprf

\section{A theory for weakly inaccessible ordinals}\label{sect:Ztheory}

Referring Theorems \ref{th:cofinalitylocal} and \ref{th:cofinalitylocalpower} let us interpret $\sfZF$ to another theory.
The base language here is $\{\in\}$.



In the following Definition \ref{df:regext}, 
$I$ is {\it intended\/} to denote  the least weakly inaccessible cardinal
though we {\it do not assume\/} the existence of weakly inaccessible cardinals
anywhere in this paper except in the {\bf Remark} after Theorem \ref{th:mainthZ}.
$\kap,\lam,\rho$ range over uncountable regular ordinals $< I$.
The predicate $P$ is intended to denote the relation
$P(\lam,x,y)$ iff $x=F^{\Sig_{1}}(\lam;I,x\cup\{\lam\})$ and 
$y=F^{\Sig_{1}}(I;I,x\cup\{\lam\})$,
and the predicate $P_{I,n}(x)$ is intended to denote the relation
$P_{I,n}(x)$ iff $x=F^{\Sig_{n}}(I;I,x)$,
where
$F^{\Sig_{n}}_{\alp}(y)=F^{\Sig_{n}}(y;I,\alp)$
denotes the Mostowski collapsing
$F^{\Sig_{n}}_{\alp}: \mbox{{\rm Hull}}_{\Sig_{n}}^{I}(\alp)\lrarw L_{\gam}$
of the $\Sig_{n}$-Skolem hull $\mbox{{\rm Hull}}_{\Sig_{n}}^{I}(\alp)$ of $\alp<I$ over $L_{I}$,
and
$F^{\Sig_{n}}_{\alp}(I):=\gam$ for $L_{\gam}=rng(F^{\Sig_{n}}_{\alp})$.


\bdf\label{df:regext}
$\mbox{{\rm T}}(I,n)$ {\rm denotes the set theory defined as follows.}
\benu
\item
{\rm Its language is} $\{\in, P,P_{I,n},Reg\}$ {\rm for a ternary predicate} $P$ {\rm  and unary predicates} $P_{I,n}$
{\rm and} $Reg$.

\item
{\rm Its axioms are obtained from those of} $\mbox{{\rm KP}}\ome$ 
{\rm
 in the expanded language
\footnote{
This means that the predicates $P,P_{I,n},Reg$ do not occur in 
$\Del_{0}$-formulae
for $\Del_{0}$-Separation and 
$\Del_{0}$-Collection.
Moreover
$P,P_{I,n}, Reg$ do not occur in Foundation axiom schema.
},
the axiom of constructibility}
$V=L$
{\rm together with the axiom schema saying that if} $Reg(\kap)$
{\rm then} $\kap$ 
 {\rm is an uncountable regular ordinal, cf. (\ref{eq:Z2}) and (\ref{eq:Z1}),
 and if} $P(\kap,x,y)$ {\rm then}  $x$ 
 {\rm is a critical point of the} $\Sig_{1}${\rm -elementary embedding from}
$L_{y}\cong \mbox{{\rm Hull}}^{I}_{\Sig_{1}}(x\cup\{\kap\})$ {\rm to the universe}
$L_{I}${\rm , cf. (\ref{eq:Z1}), and if}
$P_{I,n}(x)$ {\rm then} $x$ {\rm is a critical point of the} $\Sig_{n}${\rm -elementary embedding from}
$L_{x}\cong \mbox{{\rm Hull}}^{I}_{\Sig_{n}}(x)$ {\rm to the universe}
$L_{I}${\rm , cf.(\ref{eq:Z4}):}
{\rm for a formula} $\vphi$ {\rm and an ordinal} $\alp$,
$\vphi^{\alp}$ {\rm denotes the result of restricting every unbounded quantifier}
$\exi z,\fal z$ {\rm in} $\vphi$ {\rm to} $\exi z\in L_{\alp}, \fal z\in L_{\alp}$.

 \benu
 
 \item
 $x\in Ord$ {\rm is a} $\Del_{0}${\rm -formula saying that `}$x$ {\rm is an ordinal'.}
 \beqnarr
 &&
 (Reg(\kap) \to \ome<\kap\in Ord)
 \nonumber
\\
& \land &
(P(\kap,x,y) \to 
\{x,y\}\subset Ord \land  x<y<\kap \land Reg(\kap))
\label{eq:Z0}
\\
& \land &
(P_{I,n}(x)\to x\in Ord)
\nonumber
\eeqnarr

 \item
\beqn\label{eq:Z1}
P(\kap,x,y) \to a\in L_{x}  \to \vphi[\kap,a] \to \vphi^{y}[x,a]
\eeqn
{\rm for any} $\Sig_{1}${\rm -formula} $\vphi$ {\rm in the language} $\{\in\}$.

.

\item
\beqn\label{eq:Z2}
Reg(\kap) \to a\in Ord\cap\kap \to \exi x, y\in Ord\cap\kap[a<x\land P(\kap,x,y)]
\eeqn

\item
\beqn\label{eq:Z3}
\fal x\in Ord \exi y[y>x\land Reg(y)]
\eeqn

\item
\beqn\label{eq:Z4}
P_{I,n}(x) \to a\in L_{x} \to \vphi[a] \to \vphi^{x}[a]
\eeqn
{\rm for any} $\Sig_{n}${\rm -formula} $\vphi$ {\rm in the language} $\{\in\}$.

\item
\beqn\label{eq:Z5}
a\in Ord \to \exi x\in Ord[a<x\land P_{I,n}(x)]
\eeqn


 \eenu
 
\eenu
\edf


Let ${\sf ZFL}_{n}$ denote the subtheory of $\sfZF+(V=L)$ obtained by restricting
 {\rm Separation}  and {\rm Collection} to $\Sig_{n}$-{\rm Separation} and 
 $\Sig_{n}$-{\rm Collection}, resp.

\blem\label{lem:regularset}
$\mbox{{\rm T}}(I):=\bigcup_{n\in\ome}\mbox{{\rm T}}(I,n)$ is 
a 
conservative extension of
Zermelo-Fraenkel set theory $\sfZF+(V=L)$ with the axiom of constructibility.

Moreover for each $n\geq 1$,
$\mbox{{\rm T}}(I,n)$ is 
a 
conservative extension of
${\sf ZFL}_{n}$.
\elem
\bprf
Let $n\geq 1$.
First consider the axioms of ${\sf ZFL}_{n}$ in $T(I,n)$.
By (\ref{eq:Z4}), $\mbox{{\rm T}}(I,n)$ proves the reflection principle for $\Sig_{n}$ 
$\vphi$
\beqn\label{eq:Z6}
P_{I,n}(x) \to a\in L_{x}\to (\vphi[a]\lrarw \vphi^{x}[a])
\eeqn
Let $\vphi$ be a $\Sig_{n}$-formula, and $\alp$
 an ordinal such that $\{b,c\}\subset L_{\alp}$.
Pick an $x$ with $\alp<x\land P_{I,n}(x)$ by (\ref{eq:Z5}).
Then by (\ref{eq:Z6}) 
$\{a\in b:\vphi[a,c]\}=\{a\in b:\vphi^{x}[a,c]\}$.
This shows in $T(I,n)$, $\Sig_{n}$-{\rm Separation} from $\Del_{0}$-{\rm Separation}.
Likewise we see that $\mbox{{\rm T}}(I,n)$ proves
$\Sig_{n}$-{\rm Collection}.

Second consider the Power set axiom in $\mbox{{\rm T}}(I,n)$.
We show that the power set $\mathcal{P}(b)=\{x:x\subset b\}$ exists as a set.
Let $b\in L_{\alp}$ with a multiplicative principal number $\alp\geq\ome$. 
Pick a regular ordinal $\kap>\alp$ by (\ref{eq:Z3}).
From Theorem \ref{th:cofinalitylocal} 
we see that ${}^{\alp}\kap\subset L_{\kap}$.
Let $G:Ord\to L$ be the G\"odel's surjective map, which is $\Del_{1}$.
We have $G"\alp=L_{\alp}$ for the multiplicative principal number $\alp$.
Pick an ordinal $\bet<\alp$ such that $G(\bet)=b$.
Then
${}^{\bet}2\subset {}^{\alp}\kap\subset L_{\kap}$, i.e.,
${}^{\bet}2=\{x\in L_{\kap}: x\in {}^{\bet}2\}$, 
and hence by
$\Del_{0}$-{\rm Separation} ${}^{\bet}2$ exists as a set.
On the other hand we have
$c\in b=G(\bet) \to \exi \gam<\bet(G(\gam)=c)$ and
$\gam<\bet\to G(\gam)\in G(\bet)$.
Let $S: {}^{\bet}2\to\mathcal{P}(b)$ be the surjection
defined by
$x\in S(f)$ iff $\exi\gam<\bet(G(\gam)=x\land f(\gam)=1)$
for $f\in{}^{\bet}2$ and $x\in b$.
Pick a set $c$ such that $S"({}^{\bet}2)\subset c$ by $\Del_{0}$-Collection.
Then
$\{x:x\subset b\}=\{S(f)\in c: f\in{}^{\bet}2\}$ is a set by $\Del_{0}$-{\rm Separation}.

Hence we have shown that ${\sf ZFL}_{n}$ is contained in $\mbox{{\rm T}}(I,n)$.

Next we show that $\mbox{{\rm T}}(I,n)$ is interpretable in ${\sf ZFL}_{n}$.
Interpret the predicates
$Reg(\kap)\lrarw \ome<\kap\in Ord\land \fal \alp<\kap\fal f\in{}^{\alp}\kap[ \sup_{x<\alp}f(x)<\kap]$
and
$P(\kap,x,y)\lrarw Reg(\kap)\land \{x,y\}\subset Ord \land 
(\mbox{{\rm Hull}}_{\Sig_{1}}^{I}(x\cup\{\kap\})\cap \kap\subset x)
\land
(y=\sup\{F(a):a\in\mbox{{\rm Hull}}_{\Sig_{1}}^{I}(x\cup\{\kap\})\})$
for the Mostowski collapsing function
$F(a)=\{F(b): b\in \mbox{{\rm Hull}}_{\Sig_{1}}^{I}(x\cup\{\kap\})\cap a\}$
and the universe $L_{I}=L$.
Moreover for the predicate $P_{I,n}$,
$P_{I,n}(x)\lrarw x\in Ord \land (\mbox{{\rm Hull}}_{\Sig_{n}}^{I}(x)\cap Ord \subset x)$.

We see from Theorem \ref{th:cofinalitylocal} that
the interpreted (\ref{eq:Z0}), (\ref{eq:Z1}) and (\ref{eq:Z2})
are provable in ${\sf ZFL}_{1}$.
Moreover the unboundedness of the regular ordinals, (\ref{eq:Z3})
is provable in ${\sf ZFL}_{1}$ using the Power set axiom and $\Sig_{1}$-Separation.

It remains to show the interpreted (\ref{eq:Z4}) and (\ref{eq:Z5}) in ${\sf ZFL}_{n}$.
It suffices to show that given an ordinal $\alp$,
there exists an ordinal $x>\alp$
such that $\mbox{Hull}_{\Sig_{n}}^{I}(x)\cap Ord\subset x$.
Pick a regular ordinal $\kap>\alp$.
Again as in the proof of Proposition \ref{prp:closedpnt}.\ref{prp:closedpnt.1}, 
define recursively ordinals $\{x_{n}\}_{n}$ as follows.
$x_{0}=\alp+1$, and $x_{n+1}$ is defined to be the least ordinal $x_{n+1}$ 
such that
$\mbox{{\rm Hull}}^{I}_{\Sig_{n}}(x_{n})\cap Ord\subset x_{n+1}$.
We show inductively that such an ordinal $x_{n}$ exists, and $x_{n}<\kap$.
Then $x=\sup_{n}x_{n}\leq\kap$ is a desired one.

It suffices to show that for any $\alp<\kap$ there exists a $\bet<\kap$
such that $\mbox{Hull}_{\Sig_{n}}^{I}(\alp)\cap Ord\subset \bet$.
By Proposition \ref{prp:rhHull}
let $h^{I}_{\Sig_{n}}$ be be the $\Del_{n}$-surjection from the $\Sig_{n}$-subset
$dom(h^{I}_{\Sig_{n}})$ of $\ome\times\alp$ to
$\mbox{Hull}_{\Sig_{n}}^{I}(\alp)$, which is 
a $\Sig_{n}$-class.
From $\Sig_{n}$-Separation 
we see that $dom(h^{I}_{\Sig_{n}})$ is a set.
Hence by $\Sig_{n}$-Collection,
$\mbox{Hull}_{\Sig_{n}}^{I}(\alp)=rng(h^{I}_{\Sig_{n}})$ is a set.
Therefore the ordinal $\sup(\mbox{Hull}_{\Sig_{n}}^{I}(\alp)\cap Ord)$ exists
in the universe.
On the other hand we have for the subset $dom(h^{I}_{\Sig_{n}})$ of $\ome\times\alp$,
$dom(h^{I}_{\Sig_{n}})\in L_{\kap}$ 
by Theorem \ref{th:cofinalitylocalpower}.
Hence $\kap\leq \sup(\mbox{Hull}_{\Sig_{n}}^{I}(\alp)\cap Ord)$ would yield
a cofinal map from $\alp$ to $\kap$,
which is a subset of the set $h^{I}_{\Sig_{n}}$ in the universe.
This contradicts the regularity of $\kap$.
Therefore $\sup(\mbox{Hull}_{\Sig_{n}}^{I}(\alp)\cap Ord)<\kap$.
\eprf


\section{Ordinals for inaccessibles}\label{sect:ordinalnotation}\label{sec:ordinalinacc}
For our proof-theoretic analysis of $\sfZF+(V=L)$, we need to talk about `ordinals' 
less than the next epsilon number
to the order type of the class of ordinals inside $\sfZF+(V=L)$.
Let us define simultaneously
a $\Del_{1}$-class $Code^{\veps}$, its $\Del_{1}$-subclass $Ord^{\veps}$ and 
a $\Del_{1}$-relation $\in^{\veps}$ inside 
Kripke-Platek set theory with the axiom of infinity, $\mbox{{\rm KP}}\ome$.
On the class $Ord^{\veps}$, $\in^{\veps}$ is denoted by $<^{\veps}$, and
$x\leq^{\veps}y:\Lrarw(x=y)\lor(x<^{\veps}y)$.
$Ord^{\veps}$ is the class of codes of ordinals less than the next epsilon number
to the order type of the class of ordinals,
$Code^{\veps}$ is the union of $Ord^{\veps}$ and codes of sets in the universe.
$\in^{\veps}$ is the membership relation on codes.
Moreover we need to define two $\Del_{1}$-operations,
addition $x\oplus y$ and exponentiation $\tilde{\ome}^{x}$
on codes in $Ord^{\veps}$, which enjoy some algebraic facts demonstrably in 
$\mbox{{\rm KP}}\ome$, a fortiori in $\sfZF+(V=L)$.

Such construction is familiar in proof theory, cf. \cite{Schuette, Pohlers}.
For example in the first-order arithmetic ${\sf PA}$ we can talk abut ordinals less than $\veps_{0}$
although the order type of `ordinals' in the standard model $\mathbb{N}$ of ${\sf PA}$
is $\ome<\veps_{0}$.
We reproduce it here for completeness.
\\

For ordered pairs $\la x_{0},\ldots,x_{n}\ra$,
let $(\la x_{0},\ldots,x_{n}\ra)_{i}=x_{i}\, (i\leq n)$ and $lh(x)=n+1$.
The following Definitions \ref{df:veps} and \ref{df:oplus} are made in $\mbox{{\rm KP}}\ome$.

\bdf\label{df:veps}
\benu
\item
$Ord^{\veps}\subset Code^{\veps}$.

\item
{\rm For each set} $x$, $\la 0,x\ra\in Code^{\veps}$.
{\rm If} $x$ {\rm is an ordinal, then} $\la 0,x\ra\in Ord^{\veps}$.

\item
$\la 1,0\ra\in Ord^{\veps}$.

\item
{\rm Let} $x\in Ord^{\veps}$ {\rm be a code such that} $(x)_{0}=2,3$.
{\rm Then} $\la 2,x\ra\in Ord^{\veps}$.

\item
{\rm Let} $x_{0},\ldots,x_{m}\in Ord^{\veps}\, (m>0)$ {\rm be codes such that}
$x_{m}\leq^{\veps}\cdots\leq^{\veps}x_{0}$ {\rm and}
$[(x_{m})_{0}=0] \Rarw (x_{m})_{1}\neq 0\spand \fal i<m[(x_{i})_{0}=1,2]$.
{\rm Then}
$\la 3, x_{0},\ldots,x_{m}\ra\in Ord^{\veps}$.

\item
{\rm For} $\la 0,x\ra,\la 0,y\ra\in Code^{\veps}$,
$\la 0,x\ra\in^{\veps}\la 0,y\ra\Lrarw x\in y$.

\item
{\rm If} $x$ {\rm is an ordinal, then} 
$\la 0,x\ra<^{\veps}y$ {\rm for any} $y\in Ord^{\veps}$ {\rm with}
$(y)_{0}\neq 0$.

\item
{\rm If} $x\in Ord^{\veps}$ {\rm and} $(x)_{0}=2,3$, {\rm then} $\la 1,0\ra<^{\veps}x$.

\item
{\rm If} $x,y\in Ord^{\veps}$ {\rm and} $(x)_{0}=(y)_{0}=2$, {\rm then}
$x<^{\veps}y\Lrarw (x)_{1}<^{\veps}(y)_{1}$.

\item
{\rm If} $x,y\in Ord^{\veps}$ {\rm and} $(x)_{0}=2$, $(y)_{0}=3$, {\rm then}
$x<^{\veps}y\Lrarw x\leq^{\veps}(y)_{1}$, {\rm and}
$y<^{\veps}x\Lrarw (y)_{1}<^{\veps}x$.

\item
{\rm If} $x,y\in Ord^{\veps}$ {\rm and} $(x)_{0}=(y)_{0}=3$, {\rm then}
$x<^{\veps}y$ {\rm iff either} $lh(x)<lh(y)\spand \fal i<lh(x)[(x)_{i}=(y)_{i}]$
{\rm or} $\exi i<\min\{lh(x),lh(y)\}[(x)_{i}<^{\veps}(y)_{i}\spand \fal j<i\{(x)_{j}=(y)_{j}\}]$.
\eenu
\edf

\bdf\label{df:oplus}
\benu
\item
{\rm For} $x\in Ord^{\veps}$ {\rm and} $n<\ome$,
$\tilde{\ome}^{x}\cdot n\in Ord^{\veps}$  
{\rm is defined as follows.}
 \benu
  \item
 {\rm If} $(x)_{0}=0$, {\rm then} $\tilde{\ome}^{x}\cdot n=\la 0,\ome^{(x)_{1}}\cdot n\ra$.
 
 \item
 {\rm If} $(x)_{0}=1$, {\rm then} $\tilde{\ome}^{x}\cdot 0=\la 0,0\ra$, 
 $\tilde{\ome}^{x}\cdot 1=x$
 {\rm and for} $n>1$,
 $\tilde{\ome}^{x}\cdot n=\la 3,\la 1,0\ra,\ldots,\la 1,0\ra\ra$ {\rm with} $n${\rm 's} $\la 1,0\ra$.
 
 \item
 {\rm If} $(x)_{0}=2,3$, {\rm then} $\tilde{\ome}^{x}\cdot 0=\la 0,0\ra$,
 $\tilde{\ome}^{x}\cdot 1=\la 2,x\ra$, {\rm and for} $n>1$
 $\tilde{\ome}^{x}\cdot n=\la 3,\la 2,x\ra,\ldots,\la 2,x\ra\ra$ {\rm with} $n${\rm 's} $\la 2,x\ra$.
 \eenu

\item
$\tilde{\ome}^{x}:=\tilde{\ome}^{x}\cdot 1$.

\item
{\rm For} $x,y\in Ord^{\veps}$, $x\oplus y\in Ord^{\veps}$ {\rm is defined.}
 \benu
 \item
 {\rm If} $(x)_{0}=(y)_{0}$, {\rm then}
 $x\oplus y=\la 0, (x)_{1}+(y)_{1}\ra$.
 
 \item
  {\rm Consider the case when} $(x)_{0}=0$ {\rm and} $(y)_{0}\neq 0$.
  
   {\rm Then} $x\oplus y=y$.
   
   {\rm If} $(y)_{0}=1,2$, {\rm then} $y\oplus x=\la 3,y,x\ra$.
   
   {\rm Consider the case when} $(y)_{0}=3$ {\rm and let} $lh(y)=n+1$.
    {\rm If} $((y)_{n})_{0}=0$, {\rm then} 
    $y\oplus x=\la 3,(y)_{1},\ldots,(y)_{n\dot{-}1},(y)_{n}\oplus x\ra$.
    {\rm Otherwise} $y\oplus x=\la 3,(y)_{1},\ldots,(y)_{n}, x\ra$.
    
  \item
  {\rm Consider the case when} $\{(x)_{0},(y)_{0}\}\subset\{1,2\}$.
 {\rm If} $x<^{\veps}y$, {\rm then} $x\oplus y=y$.
 {\rm Otherwise} $x\oplus y=\la 3,x,y\ra$.
 
 \item
  {\rm Consider the case when} $(x)_{0}=1,2$ {\rm and} $(y)_{0}=3$ {\rm with} $lh(y)=n+1$.
  
  {\rm If} $x<^{\veps}(y)_{1}$, {\rm then} $x\oplus y=y$.
 {\rm Otherwise} $x\oplus y=\la 3,x,(y)_{1},\ldots,(y)_{n}\ra$.
 
 {\rm If there exists a positive integer such that}
 $k\leq n$ {\rm and} $x\leq^{\veps}(y)_{k}$.
 {\rm Then}
 $y\oplus x=\la 3, (y)_{1},\ldots,(y)_{k},x\ra$ {\rm for the largest such integer} $k$.
 {\rm Otherwise} $y\oplus x=x$.
 
 \item
 {\rm Consider the case when} $(x)_{0}=(y)_{0}=3$ {\rm with} $lh(x)=m+1$ {\rm and} $lh(y)=n+1$.
  
 {\rm If} $x<^{\veps}(y)_{1}$, {\rm then} $x\oplus y=y$.
 {\rm Otherwise let}
 $k\leq m$ {\rm denote the largest positive integer such that} $(y)_{1}\leq^{\veps}(x)_{k}$.
 {\rm Then}
 $x\oplus y=\la 3,(x)_{1},\ldots,(x)_{k},(y)_{1},\ldots,(y)_{n}\ra$.
 \eenu
\eenu
\edf

$<^{\veps}$ is seen to be a canonical ordering as stated in
the following Proposition \ref{prp:canonical}.

\bprp\label{prp:canonical}
\benu
\item\label{prp:canonical0}
$\mbox{{\sf KP}}\ome$ proves the fact that $<^{\veps}$ is a linear ordering.

\item\label{prp:canonical1}
$\mbox{{\sf KP}}\ome$ proves the facts
$x\oplus\la 0,0\ra=\la 0,0\ra\oplus x=x$, 
$y<^{\veps}x\oplus\la 0,1\ra\Lrarw y\leq^{\veps}x$
and
$x<^{\veps}y\oplus\tilde{\ome}^{z}\spand z\neq\la 0,0\ra
\Rarw \exi u<^{\veps}z\exi n<\ome(x<^{\veps}y\oplus\tilde{\ome}^{u}\cdot n)$
for $x,y,z\in Ord^{\veps}$.

\item\label{prp:canonical2}
For any formula $\vphi$
and each $n<\ome$,
\beqn\label{eq:trindveps}
\mbox{{\sf KP}}\ome\vdash\fal x\in Ord^{\veps}(\fal y<^{\veps}x\,\vphi(y)\to\vphi(x)) \to 
\fal x<^{\veps}\ome_{n}(I+1)\vphi(x)
\eeqn
where
$\ome_{n}(I+1):=\tilde{\ome}_{n}( \la 3,\la 1,0\ra, \la 0,1\ra\ra):=
\la 2,\la 2,\cdots\la 2, \la 3,\la 1,0\ra, \la 0,1\ra\ra\ra\cdots\ra\ra$,
in which
$2$ occurs $n$-times.
\eenu
\eprp
\bprf
\\
\ref{prp:canonical}.\ref{prp:canonical0} and \ref{prp:canonical}.\ref{prp:canonical1}.
It is easy to see that $<^{\veps}$ is a linear ordering, and
$\oplus, \tilde{\ome}^{x}$ enjoys the proposition.
Obviously $\la 0,0\ra$ is the least element in $Ord^{\veps}$.
\\

\noindent
\ref{prp:canonical}.\ref{prp:canonical2} is seen as usual, 
but we give a proof of it for readers' convenience.

By metainduction on $n<\ome$, we show (\ref{eq:trindveps}).
First by the Foundation axiom schema we have for any formula $\vphi$,
$\fal x\in Ord(\fal y<x\,\vphi(\la 0,y\ra)\to\vphi(\la 0,x\ra)) \to 
\fal x\in Ord\, \vphi(\la 0,x\ra)$.
This yields the case $n=0$.

Let $Prg[\vphi]:\Lrarw \fal x\in Ord^{\veps}(\fal y<^{\veps}x\,\vphi(y)\to\vphi(x))$,
and suppose $Prg[\vphi]\to \fal x<^{\veps}\ome_{n}(I+1)\vphi(x)$ for any formula $\vphi$.

Given a formula $\vphi(x)$, let
${\sf j}[\vphi](a):\Lrarw 
\fal x\in Ord^{\veps}(\fal y<^{\veps}x\,\vphi(y)\to \fal y<^{\veps}x\oplus \tilde{\ome}^{a}\,\vphi(y))$.
Then from Proposition \ref{prp:canonical}.\ref{prp:canonical1}
we see that $Prg[\vphi]\to Prg[{\sf j}[\vphi]]$.
Assuming $Prg[\vphi]$, we obtain $Prg[{\sf j}[\vphi]]$.
IH yields
$\fal x<^{\veps}\ome_{n}(I+1){\sf j}[\vphi](x)$, and hence
${\sf j}[\vphi](\ome_{n}(I+1))$.
Therefore by $\la 0,0\ra\oplus z=z$ we conclude
$ \fal y<^{\veps}\tilde{\ome}^{\ome_{n}(I+1)}\,\vphi(y)$ 
for $\tilde{\ome}^{\ome_{n}(I+1)}=\ome_{n+1}(I+1)$.
\eprf
\\

In what follows of this section $n\geq 1$ denotes a fixed positive integer, 
and we work in $\sfZF+(V=L)$.

Let
\[
I:=\la 1,0\ra, \: \ome_{n}(I+1):=\tilde{\ome}_{n}( \la 3,\la 1,0\ra, \la 0,1\ra\ra), 
\mbox{ and } L_{I}:=\{\la 0,x\ra: x\in L\}
\]
and for codes $X,Y\in Code^{\veps}$ 
\[
X\subset^{\veps}Y:\Lrarw \fal x\in^{\varepsilon} X(x\in^{\varepsilon}Y)
.\]

For simplicity let us identify the code $x\in Code^{\veps}$ with
the `set' coded by $x$,
and $\in^{\veps}$ [$<^{\veps}$] is denoted by $\in$ [$<$], resp. when no confusion likely occurs.
For example, the code $\la 0,x\ra$ is identified with the set $\{\la 0,y\ra: y\in x\}$ of codes.

$R:=\{\la 0,\kap\ra: \kap \mbox{ is a uncountable regular ordinal}\}$, while
$R^{+}:=R\cup\{I\}$.
$\kap,\lam,\rho$ denote elements of $R$.

Define simultaneously 
 the classes $\calh_{\alp,n}(X)\subset L_{I}\cup\{x\in Ord^{\veps}: x<^{\veps}\ome_{n+1}(I+1)\}$
and the ordinals $\Psi_{\kap,n} \alp\,(\kap\in R)$ and $\Psi_{I,n}\alp$ 
for $\alp<^{\veps}\ome_{n+1}(I+1)$ and {\it sets\/} $X\subset L_{I}$ as follows.
We see that $\calh_{\alp,n}(X)$ and $\Psi_{\kap,n} \alp$ are (first-order) definable as a fixed point in $\sfZF$, cf. Proposition \ref{prp:definability}.

\bdf\label{df:Cpsiregularsm}

$\calh_{\alp,n}(X)$ {\rm is the} 
{\rm Skolem hull of} $\{\langle 0,0\rangle,I\}\cup X$ {\rm under the functions} 
$\oplus,
 \alp\mapsto\tilde{\ome}^{\alp}<^{\veps}\ome_{n+1}(I+1), 
 \Psi_{I,n}\restrict\alp, \Psi_{\kap,n}\restrict \alp\,(\kap\in R)${\rm , the}
 $\Sig_{n}${\rm -definability,}
{\rm and the Mostowski collapsing functions}
$(x,\kap,d)\mapsto F^{\Sig_{1}}_{x\cup\{\kap\}}(d)\,
(\kap\in R, \mbox{{\rm Hull}}^{I}_{\Sig_{1}}(x\cup\{\kap\})\cap\kap\subset x)$
{\rm and}
$(x,d)\mapsto F^{\Sig_{n}}_{x}(d)\,
(\mbox{{\rm Hull}}^{I}_{\Sig_{n}}(x)\cap I\subset x)$.

{\rm For a later reference let us define stages}
$\calh_{\alp,n}^{m}(X)\,(m\in\ome)$ {\rm of the inductive definition.}

\benu
\item
$\calh_{\alp,n}^{0}(X)=\{\langle 0,0\rangle,I\}\cup X$.

\item
 $x, y \in \calh_{\alp,n}^{m}(X)\cap \ome_{n+1}(I+1) \Rarw x\oplus y\in \calh_{\alp,n}^{m+1}(X)$.
 \\
 $x \in \calh_{\alp,n}^{m}(X)\cap \ome_{n}(I+1) \Rarw \tilde{\ome}^{x}\in \calh_{\alp,n}^{m+1}(X)$.


\item
$\gam\in \calh_{\alp,n}^{m}(X)\cap\alp
\Rarw 
\Psi_{I,n}\gam\in\calh_{\alp,n}^{m+1}(X)
$.

\item
$
\kap\in\calh_{\alp,n}^{m}(X)\cap R
\spand
\gam\in \calh_{\alp,n}^{m}(X)\cap\alp
\Rarw 
\Psi_{\kap,n}\gam\in\calh_{\alp,n}^{m+1}(X)
$.

\item
$
\mbox{{\rm Hull}}^{I}_{\Sig_{n}}(\calh_{\alp,n}(X)\cap L_{I})\cap Code^{\veps} 
\subset \calh_{\alp,n}(X)
$.

{\rm Namely for any} $\Sig_{n}${\rm -formula} $\vphi[x,\vec{y}]$ {\rm in the language} $\{\in\}$
{\rm and parameters} $\vec{a}\subset \calh_{\alp,n}^{m}(X)\cap L_{I}${\rm , if} 
$b\in L_{I}$,
$(L_{I},\in^{\veps})\models\vphi[b,\vec{a}]$ {\rm and} 
$(L_{I},\in^{\veps})\models\exi!x\,\vphi[x,\vec{a}]${\rm , then} $b\in\calh_{\alp,n}^{m+1}(X)$.


\item
{\rm If} $\kap\in\calh_{\alp,n}^{m}(X)\cap R$, 
$x\in\calh_{\alp,n}^{m}(X)\cap \kap$ {\rm with} 
$\mbox{{\rm Hull}}^{I}_{\Sig_{1}}(x\cup\{\kap\})\cap\kap\subset^{\veps} x$ {\rm and}
$(\kap=\ome_{c+1}\Rarw \ome_{c}<x)$, 
{\rm and}
$d\in (\mbox{{\rm Hull}}^{I}_{\Sig_{1}}(x\cup\{\kap\})\cup\{I\})\cap\calh_{\alp,n}^{m}(X)$,
{\rm then}
$F^{\Sig_{1}}_{x\cup\{\kap\}}(d)\in\calh_{\alp,n}^{m+1}(X)$.

\item
{\rm If} 
$x\in\calh_{\alp,n}^{m}(X)\cap I$ {\rm with} 
$\mbox{{\rm Hull}}^{I}_{\Sig_{n}}(x)\cap I\subset^{\veps}  x$,
{\rm and}
$d\in (\mbox{{\rm Hull}}^{I}_{\Sig_{n}}(x)\cup\{I\})\cap\calh_{\alp,n}^{m}(X)$,
{\rm then}
$F^{\Sig_{n}}_{x}(d)\in\calh_{\alp,n}^{m+1}(X)$.

\item
$
\calh_{\alp,n}(X):=\bigcup\{\calh_{\alp,n}^{m}(X):m\in\ome\}
$.

\eenu


{\rm For} $\kap\in R^{+}$
\[
\Psi_{\kap,n}\alp:=
\min_{<^{\veps}}(\{\kap\}\cup\{\bet<^{\veps}\kap : \kap\in \calh_{\alp,n}(\bet)\spand \calh_{\alp,n}(\bet)\cap \kap \subset^{\veps} \bet\})
.\]


\edf

The ordinal $\Psi_{\kap,n}\alp$ is well defined and $\Psi_{\kap,n}\alp\leq^{\veps} \kap$
for any uncountable regular $\kap\leq^{\veps}  I$ 
since $\kap\in\calh_{\alp,n}(\kap)$ 
by Proposition \ref{prp:lem4.5}.\ref{prp:lem4.5.1} below.


\bprp\label{prp:clshull}
\benu
\item\label{prp:clshull.0}
$\calh_{\alp,n}(X)$ is closed under $\Sig_{n}$-definability:
$
\vec{a}\subset\calh_{\alp,n}(X)\cap L_{I}  \Rarw \mbox{{\rm Hull}}^{I}_{\Sig_{n}}(\vec{a})\cap L_{I}\subset \calh_{\alp,n}(X)
$.

\item\label{prp:clshull.2}
For $\kap\in R$,
$
\mbox{{\rm Hull}}^{I}_{\Sig_{1}}(\Psi_{\kap,n}\alp\cup\{\kap\})\cap\kap\cap Code^{\veps} =\Psi_{\kap,n}\alp
$.
Namely $\Psi_{\kap,n}\alp\in Cr^{I}_{\Sig_{1}}(\{\kap\})$.

\item\label{prp:clshull.1}
$\calh_{\alp,n}(X)$ is closed under the Veblen function $\tilde{\vphi}$ on $I$,
$x, y \in \calh_{\alp,n}(X)\cap I \Rarw \tilde{\vphi} xy\in \calh_{\alp,n}(X)$,
where
$\tilde{\vphi}(\la 0,x\ra)(\la 0, y\ra)=\la 0,\vphi xy\ra$ for ordinals $x,y$.

\item\label{prp:clshull.3}
If
$\kap\in\calh_{\alp,n}(X)\cap R$, $x\in\calh_{\alp,n}(X)\cap\kap$, $\mbox{{\rm Hull}}^{I}_{\Sig_{1}}(x\cup\{\kap\})\cap\kap\subset^{\veps} x$,
$(\kap=\ome_{c+1}\Rarw \ome_{c}<x)$ and
$\del\in(\mbox{{\rm Hull}}^{I}_{\Sig_{1}}(x\cup\{\kap\})\cup\{I\})\cap \calh_{\alp,n}(X)$,
then
$F^{\Sig_{1}}_{x\cup\{\kap\}}(\del)\in\calh_{\alp,n}(X)$.

\item\label{prp:clshull.3n}
If
$x\in\calh_{\alp,n}(X)\cap I$, $\mbox{{\rm Hull}}^{I}_{\Sig_{n}}(x)\cap I\subset^{\veps} x$ and
$\del\in(\mbox{{\rm Hull}}^{I}_{\Sig_{n}}(x)\cup\{I\})\cap \calh_{\alp,n}(X)$,
then
$F^{\Sig_{n}}_{x}(\del)\in\calh_{\alp,n}(X)$.

\item\label{prp:clshull.4}
Assume $n\geq 2$.
$\gam\in \calh_{\alp,n}(X)\cap I
\Lrarw
\ome_{\gam}\in \calh_{\alp,n}(X)\cap I
$ 
for $\ome_{\alp}=\aleph_{\alp}$.

Moreover 
$\gam\in \calh_{\alp,n}(X)\cap I
\Rarw
\gam^{+}=\min\{\lam\in R: \gam<\lam\}\in \calh_{\alp,n}(X)\cap I
$ .
\eenu
\eprp
\bprf
\ref{prp:clshull}.\ref{prp:clshull.2}.
By the definition of $\calh_{\alp,n}(X)$, we have
\[
\mbox{{\rm Hull}}^{I}_{\Sig_{1}}(\Psi_{\kap,n}\alp\cup\{\kap\})\cap\kap
\subset^{\veps} 
\calh_{\alp,n}(\Psi_{\kap,n}\alp)\cap\kap\subset^{\veps} 
\Psi_{\kap,n}\alp
\subset
\mbox{{\rm Hull}}^{I}_{\Sig_{1}}(\Psi_{\kap,n}\alp\cup\{\kap\})\cap\kap
.\]
\ref{prp:clshull}.\ref{prp:clshull.1}.
This is seen from the $\Sig_{1}$-definability of the Veblen function $\vphi$.
\\
\ref{prp:clshull}.\ref{prp:clshull.4}.
From Corollary \ref{cor:cofinalreg}.\ref{cor:cofinalreg15} 
the map $I>\alp\mapsto\ome_{\alp}$ and its inverse
are $\Del_{2}$-definable.
Moreover the next regular ordinal $\gam^{+}$ is $\Del_{2}$-definable.
\eprf
\\

In the following Proposition \ref{prp:DeldefF},
 for $\kap\in R^{+}$ and $x$,
 $(\mbox{Hull}(x,\kap),F_{x,\kap})$ denotes $(\mbox{Hull}^{I}_{\Sig_{1}}(x\cup\{\kap\}), F^{\Sig_{1}}_{x\cup\{\kap\}})$ when $\kap<^{\veps} I$, and $(\mbox{Hull}^{I}_{\Sig_{n}}(x), F^{\Sig_{n}}_{x})$ when $\kap=I$.

\bprp\label{prp:DeldefF}
Suppose $n\geq 2$, $\kap,\lam\in R^{+}$, $\mbox{{\rm Hull}}(x,\kap)\cap\kap\subset^{\veps} x$,
and $\la 0,\ome_{c}\ra<^{\veps}x$ if $\kap=\la 0,\ome_{c+1}\ra$. Then
$
x<^{\veps} \Psi_{\lam,n} b \Rarw F_{x,\kap}(I)<^{\veps} \Psi_{\lam,n} b
$,
and
$
a\in\calh_{b,n}(\Psi_{\kap,n}b)\cap b \Rarw \Psi_{\kap,n}a<^{\veps} \Psi_{\kap,n}b
$.
\eprp
\bprf
Suppose 
$x<^{\veps} \Psi_{\lam,n} b$.
We show $\kap\in\calh_{b,n}(\Psi_{\lam,n}b)$.
If $\kap=I$, there is nothing to show.
If $\kap=\la 0,\ome_{c+1}\ra$, we have $\la 0,c\ra\leq^{\veps} \la 0,\ome_{c}\ra<^{\veps} x<^{\veps} \Psi_{\lam,n}b$ and
$\la 0,c\ra\in\calh_{b,n}(\Psi_{\lam,n}b)$.
By Proposition \ref{prp:clshull}.\ref{prp:clshull.4} we have 
$\kap=\la 0,\ome_{c+1}\ra\in\calh_{b,n}(\Psi_{\lam,n}b)$.
Thus $F_{x,\kap}(I)\in\calh_{b,n}(\Psi_{\lam,n}b)$.
It remains to see $y:=F_{x,\kap}(I)<^{\veps} \lam$.
%
%
%
We have a definable bijection from $x$ to $L_{y}$.
Since $x<^{\veps} \lam$, we conclude $F_{x,\kap}(I)=y<^{\veps} \lam$.
\eprf
\\

We see the following Proposition \ref{prp:lem4.5} as in \cite{Buchholz}.

\bprp\label{prp:lem4.5}
Let $n\geq 2$.
\benu
\item\label{prp:lem4.5.1}
For any $\kap\in R^{+}$,
$\kap\in\calh_{\alp,n}(\kap)$, $\kap\in\calh_{\alp,n}(\Psi_{\kap,n}\alp)$ and $\Psi_{\kap,n}\alp<\kap$.

\item\label{prp:lem4.5.2}
$\Psi_{\kap,n}\alp\not\in\{\ome_{\bet}:\bet<\ome_{\bet}\}$.

\item\label{prp:lem4.5.3}
$\ome_{a}<\Psi_{\ome_{a+1},n}\alp<\ome_{a+1}$.

\item\label{prp:lem4.5.4}
$\ome_{\Psi_{I,n}\alp}=\Psi_{I,n}\alp$.

\item\label{prp:lem4.5.5}
$\Psi_{I,n}\alp<I$

\eenu
\eprp


The following Proposition \ref{prp:definability} is easy to see.

\bprp\label{prp:definability}
Both of 
$x=\calh_{\alp,n}(X)\, (\alp<^{\veps} \ome_{n+1}(I+1), X\in L_{I})$ and
$y=\Psi_{\kap,n}\alp\,(\kap\in R^{+})$
are $\Sig_{n+1}$-predicates as fixed points in $\sfZF$. 
\eprp


\blem\label{lem:lowerbndreg}
$\fal\alp<^{\veps}\ome_{n+1}(I+1)
\fal\kap\in R^{+}\exi x<^{\veps}\kap[x=\Psi_{\kap,n}\alp]
$.
\elem
\bprf
By Proposition \ref{prp:definability}
 both 
$x=\calh_{\alp,n}(\bet)\, (\alp<^{\veps}\ome_{n+1}(I+1),\bet<^{\veps}I)$ and
$y=\Psi_{\kap,n}\alp\,(\kap\in R^{+})$ are $\Sig_{n+1}$-predicates.
We show
that 
$A(\alp) :\Lrarw
 \fal\bet<^{\veps}I\exi x[x=\calh_{\alp,n}(\bet)] \land \fal\kap\in R^{+}\exi\bet<^{\veps}\kap[\Psi_{\kap,n}\alp=\bet]$ is 
progressive along $<^{\veps}$.
Then $\fal\alp<^{\veps}\ome_{n+1}(I+1) \fal\kap\in R^{+}\exi x<^{\veps}\kap[x=\Psi_{\kap,n}\alp]$ will follow from transfinite induction up to $\ome_{n+1}(I+1)$, cf. (\ref{eq:trindveps}) in Proposition \ref{prp:canonical}.

Assume $\fal\gam<^{\veps}\alp\, A(\gam)$ as our IH.
We have $\fal x\exi h[h=\mbox{{\rm Hull}}_{\Sig_{n}}^{I}(x)]$.
We see from this, IH and 
Separation that $\fal X\exi ! Y\, D_{\alp,n}(X,Y)$,
where $D_{\alp,n}(X,Y)$ is a $\Sig_{n+1}$-predicate such that
if $D_{\alp,n}(\calh_{\alp,n}^{m}(\bet),Y)$ then $Y=\calh_{\alp,n}^{m+1}(\bet)$ 
for any $Y$.
Therefore $\fal\bet<^{\veps}I\exi x[x=\calh_{\alp,n}(\bet)=\bigcup_{m}\calh_{\alp,n}^{m}(\bet)]$.

Next as in the Proof of Theorem \ref{th:cofinalitylocal}, (\ref{eqn:cofinalitylocal3}) $\Rarw$ (\ref{eqn:cofinalitylocal0}),
define recursively codes of ordinals $\{\bet_{m}\}_{m}$ for $\kap\in R^{+}$ as follows.
$\bet_{0}=\la 0,0\ra$ if $\kap=I$ and $\bet_{0}=\la 0, a+1\ra$ if $\kap=\la 0,\ome_{a+1}\ra$, and 
$\bet_{m+1}$ is defined to be the code of the least ordinal $\bet_{m+1}\leq^{\veps}\kap$ such that
$\calh_{\alp,n}(\bet_{m})\cap\kap\subset^{\veps}\bet_{m+1}$.

We see inductively that $\bet_{m}<^{\veps}\kap$ using the regularity of $\kap$ and the facts that
$\fal \bet<^{\veps}I\exi x[x=\calh_{\alp}(\bet)]$ and 
$\fal\bet<^{\veps}\kap\exi x[x=\calh_{\alp}(\bet)
 \land card(x)<(\kap)_{1}]$ for $\kap=\la 0,\ome_{a+1}\ra$.
 For the case $\kap=I$, $card(x)<(\kap)_{1}$ should be replaced by $card(x)<\ome_{1}$.
The latter follows from the fact that 
$\fal X\exi Y[D_{\alp,n}(X,Y)\land \{card(X)<\kap \to card(Y)<\kap\}]$.

Moreover $m\mapsto\bet_{m}$ is a definable map.
Therefore $\bet=\sup_{m}\bet_{m}<^{\veps}\kap$ enjoys 
$\calh_{\alp,n}(\bet)\cap\kap\subset^{\veps}\bet$.
Also $\la 0,a\ra\in\calh_{\alp,n}(\bet)$ for $\kap=\la 0,\ome_{a+1}\ra$.
\eprf

\section{Operator controlled derivations for weakly inaccessibles}\label{sect:controlledOme}

This section relies on Buchholz' techniques in \cite{Buchholz}.

In what follows of this section $n\geq 2$ denotes a fixed positive integer.
We consider only the codes of the ordinals $<^{\veps}\ome_{n+1}(I+1)$.

For $\alp<^{\veps}I=\la 1,0\ra$, $L_{\alp}=\{\la 0,x\ra:x\in L_{(\alp)_{1}}\}$.
$L_{I}=\{\la 0,x\ra:x\in L\}=\bigcup_{\alp<^{\veps}I}L_{\alp}$
denotes the universe.
Both $(L_{I},\in^{\veps})\models A$ and `$A$ is true' are synonymous with $A$.

\subsection{An intuitionistic fixed point theory $\mbox{FiX}^{i}(\sfZFL)$}\label{subsec:intfixZFL}
To analyze the theory \sfZF+(V=L), we first embed finitary derivations of $\exi x\in L_{\ome_{1}}\,\vphi(x)$
in \sfZF+(V=L) to
infinitary derivations, called {\it operator controlled derivations\/} due to Buchholz\cite{Buchholz}.
And then cut inferences are partially eliminated, and 
$\exi x\in L_{\ome_{1}}\,\vphi(x)$ is collapsed down to 
countables $\exi x\in L_{\Psi_{\ome_{1},n}\ome_{n}(I+1)} \vphi(x)$ for an $n<\ome$.
From the resulting derivation we need to conclude that 
the end formula $\exi x\in L_{\Psi_{\ome_{1},n}\ome_{n}(I+1)} \vphi(x)$ is true in \sfZF+(V=L).

The central notion is the relation
$(\calh_{\gam}[\Tht],\kap,n)\vdash^{a}_{b}\Gam$ defined in subsection \ref{subsec:opcontderivation},
where $n$ is a fixed integer, 
$\gam,\kap,a,b$ are codes of ordinals with $a<^{\veps}\ome_{n}(I+1)$, $b<^{\veps}I\oplus\ome$ and $\kap\leq^{\veps} I$ the code of a regular ordinal,
$\Tht$ is a finite subset of $L_{I}$ and $\Gam$ a sequent, i.e., a finite set of sentences.
Usually the relation is defined by recursion on `ordinals' $a$, but such a recursion is not available in \sfZF+(V=L)
since $a$ may be larger than $I$.
Instead of the recursion, the relation is defined for each $n<\ome$, as a fixed point,
\beqn\label{eq:fixH}
H_{n}(\gam,\Tht,\kap,a,b,\Gam)\Lrarw (\calh_{\gam,n}[\Tht],\kap,n)\vdash^{a}_{b}\Gam
\eeqn
Once this is done, it suffices to have, besides \sfZF+(V=L), the axiom (\ref{eq:fixH})
and transfinite induction schema up to each ordinal$<^{\veps}\ome_{n}(I+1)$
applied to any (first-order) formulas in the language $\{\in,H_{n}\}$ 
to carry out the proofs in this section.
Now a crucial observation due to Buchholz \cite{BuchholzAFML} 
is that the underlying logic in the proofs is {\it intuitionistic\/}.
This means that the whole proof in this section is formalizable in 
an intuitionistic fixed point theory $\mbox{FiX}^{i}(\sfZFL)$ over \sfZF+(V=L).
Then as in \cite{intfix}
we see that $\mbox{FiX}^{i}(\sfZFL)$ is a conservative extension of \sfZF+(V=L).
In this way we can conclude that 
the end formula $\exi x\in L_{\Psi_{\ome_{1},n}\ome_{n}(I+1)} \vphi(x)$ is true in \sfZF+(V=L).
\\

Let $T$ be a recursive set theory in the language $\{\in,=\}$ containing $\mbox{KP}\ome$,
where by saying that 
$T$ is a recursive theory we mean that the set of G\"odel numbers of axioms in $T$
is recursive under a standard encoding of syntax.
$\mbox{KP}\ome,\sfZF,\sfZF+(V=L)$ are examples of such theories $T$.

Let us introduce an intuitionistic fixed point theory $\mbox{FiX}^{i}(T)$ over $T$.
When $T=\sfZF+(V=L)$, we write $\mbox{FiX}^{i}(\sfZFL)$ for $\mbox{FiX}^{i}(\sfZF+(V=L))$.
Fix an $X$-strictly positive formula $\mathcal{Q}(X,x)$ in the language $\{\in,=,X\}$ with an extra unary predicate symbol $X$.
In $\mathcal{Q}(X,x)$ the predicate symbol $X$ occurs only strictly positive.
This means that the predicate symbol $X$ does not occur in the antecedent $\vphi$ of implications $\vphi\to\psi$ 
nor in the scope of negations $\lnot$ in $\mathcal{Q}(X,x)$.
The language of $\mbox{FiX}^{i}(T)$ is $\{\in,=,Q\}$ with a fresh unary predicate symbol $Q$.
The axioms in $\mbox{FiX}^{i}(T)$ consist of the following:
\benu
\item
All provable sentences in $T$ (in the language $\{\in,=\}$).

\item
Induction schema for any formula $\vphi$ in $\{\in,=,Q\}$:
\beqn\label{eq:Qind}
\fal x(\fal y\in x\,\vphi(y)\to\vphi(x))\to\fal x\,\vphi(x)
\eeqn

\item
Fixed point axiom:
\[
\fal x[Q(x)\lrarw \mathcal{Q}(Q,x)]
.\]
\eenu

The underlying logic in $\mbox{FiX}^{i}(T)$ is defined to be the intuitionistic (first-order predicate) logic (with equality).

\blem\label{lem:vepsfix}
Let $<^{\veps}$ denote a $\Del_{1}$-predicate defined in section \ref{sec:ordinalinacc}.
For each $n<\ome$ and each formula $\vphi$ in $\{\in,=,Q\}$,
\[
\mbox{{\rm FiX}}^{i}(T)\vdash\fal x(\fal y<^{\veps}x\,\vphi(y)\to\vphi(x)) \to 
\fal x<^{\veps}\ome_{n}(I+1)\vphi(x)
.\]
\elem
\bprf
This is seen as in (\ref{eq:trindveps}) of Proposition \ref{prp:canonical} using (\ref{eq:Qind}).
\eprf

\bth\label{th:consvintfix}
For any recursive set theory $T\supset\mbox{KP}\ome$,
$\mbox{{\rm FiX}}^{i}(T)$ is a conservative extension of $T$.
\end{theorem}
\bprf
Since this is seen as in \cite{intfix}, our proof is sketchy.

First the finitary derivations of set-theoretic sentences $\vphi$ in $\mbox{{\rm FiX}}^{i}(T)$ are embedded to infinitary derivations of a sequent
$\tht\Rarw\vphi$ for a provable sentence $\tht$ in $T$,
where in a sequent $\Gam\Rarw A$, $\Gam\cup\{A\}$ is a set of sentences in the language
$\{\in,=,Q\}\cup V^{\veps}$ with $V^{\veps}=\{\la 0,x\ra: x \mbox{ is a set}\}$.
Each $\la 0,x\ra\in V^{\veps}$ is an individual constant denoting $x$.
Initial sequents in infinitary derivations are
\[
\Gam,Q(a)\Rarw Q(a); \:\: \Gam,\bot\Rarw A; \:\: \Gam\Rarw \top
\]
where $a\in V^{\veps}$ and 
$\bot$ [$\top$] stands ambiguously for false atomic sentences
[true atomic sentences] in the languauge $\{\in,=\}\cup V^{\veps}$, resp.

Inference rules in infinitary derivations are standard ones for propositional connectives $\lor,\land,\to$ and
the followings for $n<\ome$,
\[
\infer[(LQ)]
{\vdash^{\alp}_{n}\Gamma,Q(a) \Rightarrow C}
{\vdash^{\bet}_{n}\Gamma,Q(a), \mathcal{Q}(Q,a) \Rightarrow C}
\: ;\:
\infer[(RQ)]
{\vdash^{\alp}_{n}\Gamma\Rightarrow Q(a)}
{\vdash^{\bet}_{n}\Gamma\Rightarrow \mathcal{Q}(Q,a)}
\]
for $\bet<^{\veps}\alp$ and any $a\in V^{\veps}$.
 \[
\infer[(R\exists)]
{\vdash^{\alp}_{n}\Gamma\Rightarrow \exists x B(x)}
{\vdash^{\bet}_{n}\Gamma\Rightarrow B(a)}
 \: ;\:
 \infer[(L\forall)]
 {\vdash^{\alp}_{n}\Gamma,\forall x B(x)\Rightarrow C}
 {\vdash^{\bet}_{n}\Gamma,\forall x B(x),B(a)\Rightarrow C}
 \]
for $\bet<^{\veps}\alp$ and any $a\in V^{\veps}$.
\[
\infer[(L\exists)]
{\vdash^{\alp}_{n}\Gamma,\exists x B(x)\Rightarrow C}
{
\cdots
&
\vdash^{\bet_{b}}_{n}\Gamma,\exists x B(x), B(b)\Rightarrow C
&
\cdots (b\in V^{\veps})
}
\: ;\:
 \infer[(R\forall)]
 {\vdash^{\alp}_{n}\Gamma\Rightarrow \forall x B(x)}
 {
 \cdots
 &
\vdash^{\bet_{b}}_{n} \Gamma\Rightarrow B(b)
 &
 \cdots (b\in V^{\veps})
 }
 \]
 where $\fal b\in V^{\veps}(\bet_{b}<^{\veps}\alp)$.
  \[
  \infer[(cut)]
  {\vdash^{\alp}_{n}\Gamma\Rightarrow C}
  {
  \vdash^{\bet}_{n}\Gamma\Rightarrow A
  &
 \vdash^{\bet}_{n}\Gam,A\Rightarrow C
  }
 \: ;\:
  \infer[(Rep)]
  {\vdash^{\alp}_{n}\Gamma\Rightarrow C}{ \vdash^{\bet}_{n}\Gamma\Rightarrow C}
  \]
where $\bet<^{\veps}\alp$, and the number of occurrences of logical symbols 
$\{\lor,\land,\to,\exi,\fal\}$ in the cut formula $A$ is smaller than $n$.

From the infinitary derivation, eliminate cuts partially.
This results in an infinitary derivation of the same sequent $\tht\Rarw\vphi$
with depth$<^{\veps}\ome_{m}(I+1)$ for an $m<\ome$, which depends on the given finite derivation.
In the resulted derivation there
occur cut inferences with cut formulas $Q(x)$ for fixed points only.
Now the constraint on operator $\mathcal{Q}$ admits us to
eliminate strictly positive cut formulas quickly as in \cite{intfix}.
In this way we will get an infinitary derivation of $\tht\Rarw\vphi$
with depth$<^{\veps}\ome_{m+2}(I+1)$, in which there
occur no fixed point formulas.

By formalizing the arguments in $T$ 
we see that the end formula $\vphi$ is true in $T$.
In the formalization, note that `$\lceil A\rceil$ is a code of a sentence $A$' is a $\Del_{1}$-predicate,
and an infinitary derivation is a class $\{T(s):s\in T\}$ 
of sequents $T(s)$ together with some informations on
the last inference and ordinal depths as in \cite{Mintsfinite},
where $T$ is a tree of finite sequences of elements in $V^{\veps}$.
Thus `$T$ is an infinitary derivation' is a $\Pi_{1}$-predicate.
Moreover any infinitary derivation in the proof can be assumed to be recursive.
This means that each infinitary derivation is encoded by a code $e\in V^{\veps}$
such that for each finite sequence $s$ of elements in $V^{\veps}$,
$\{e\}(s)$ is either $\emptyset$, which tells us $s$ is not in the derivation tree $T$,
or $\{e\}(s)=T(s)$. 
Finally by induction up to ordinal depths $\alp<^{\veps}\ome_{m+2}(I+1)$,
which is available in $T$ by (\ref{eq:trindveps}) of Proposition \ref{prp:canonical}, and
a partial truth definition of sentences in a constant logical complexity,
we see that the end sequent $\tht\Rarw\vphi$ is true, and hence so is $\vphi$ in 
$T$.
\eprf
\\

In the remaining parts of this section we work in the intuitionistic fixed point theory $\mbox{FiX}^{i}(\sfZFL)$.

\subsection{Classes of formulae}

The language $\calL_{c}$ is obtained from $\{\in, P,P_{I,n},Reg\}$ 
by adding names(individual constants) $c_{a}$
of each set $a\in L_{I}$.
$c_{a}$ is identified with $a$.
A {\it term\/} in $\calL_{c}$ is either a variable or a constant in $L_{I}$.

Formulae in this language are defined in the next definition.
Formulae are assumed to be in negation normal form.

\bdf\label{df:OTfml}
\benu
\item
{\rm Let} $t_{1},\ldots,t_{m}$ {\rm be terms.}
{\rm For each} $m${\rm -ary predicate constant} $R\in\{\in,P,P_{I,n},Reg\}$ 
$R(t_{1},\ldots,t_{m})$ {\rm and} $\lnot R(t_{1},\ldots,t_{m})$ {\rm are formulae,
where} $m=1,2,3$.
{\rm These are called} literals.


\item
{\rm If} $A$ {\rm and} $B$ {\rm are formulae, then so are} $A\land B$ {\rm and} $A\lor B$.

\item
{\rm Let} $t$ {\rm be a term.}
{\rm If} $A$ {\rm is a formula and the variable} $x$ {\rm does not occur in} $t${\rm , then}
$\exi x\in t\, A$ {\rm and} $\fal x\in t\, A$ {\rm are} bounded {\rm formulae.}

\item
{\rm If} $A$ {\rm is a formula and} $x$ {\rm a variable, then}
$\exi x\, A$ {\rm and} $\fal x\, A$ {\rm are} unbounded {\rm formulae.}
{\rm Unbounded quantifiers} $\exi x,\fal x$ {\rm are denoted by} $\exi x\in L_{I},\fal x\in L_{I}${\rm , resp.}
\eenu
\edf

For formulae $A$ in $\calL_{c}$, 
$\sfqk(A)$ denotes the finite set 
of sets 
$a\in L_{I}$ 
which are bounds of `bounded' quantifiers $\exi x\in a,\fal x\in a$
occurring in $A$.
Moreover
$\sfk(A)$ denotes the set of
sets occurring in $A$.
$\sfk(A)$ is defined to include 
bounds of `bounded' quantifiers.
By definition we set $0\in\sfqk(A)$.
Thus $0\in\sfqk(A)\subset\sfk(A)\subset L_{I}$.


\bdf
\benu
\item
$\sfk(\lnot A)=\sfk(A)$ {\rm and similarly for} $\sfqk$.
\item
$\sfqk(M)=\{0\}$ {\rm for any literal} $M$.
\item
$\sfk(Q(t_{1},\ldots,t_{m}))=(\{t_{1},\ldots,t_{m}\}\cap L_{I})\cup\{0\}$
{\rm for literals} $Q(t_{1},\ldots,t_{m})$ {\rm with predicates} $Q$ {\rm in the set} 
$\{\in,P,P_{I,n},Reg\}$.
\item
$\sfk(A_{0}\lor A_{1})=\sfk(A_{0})\cup\sfk(A_{1})$ {\rm and similarly for} $\sfqk$.
\item
{\rm For unbounded quantifiers},
$\sfk(\exi x\, A(x))=\sfk(A(x))$ {\rm and similarly for} $\sfqk$.
\item
{\rm For bounded quantifiers with} $a\in L_{I}$,
$\sfk(\exi x\in a\, A(x))=\{a\}\cup\sfk(A(x))$ {\rm and similarly for} $\sfqk$.
\item
{\rm For variables} $y$,
$\sfk(\exi x\in y\, A(x))=\sfk(A(x))$ {\rm and similarly for} $\sfqk$.
\item
{\rm For sets} $\Gam$ {\rm of formulae}
$\sfk(\Gam):=\bigcup\{\sfk(A):A\in\Gam\}$.
\eenu
\edf

For example
$\sfqk(\exi x\in a\, A(x))=\{a\}\cup\sfqk(A(x))$
if $a\in L_{I}$.

\bdf
{\rm For} $a\in L_{I}\cup\{L_{I}\}${\rm ,} $\rk_{L}(a)$ {\rm denotes the} $L$-rank {\rm of} $a$.
\[
\rk_{L}(a):=
\left\{
\begin{array}{ll}
\min\{\alp\in Ord: a\in L_{\alp+1}\} & a\in L_{I}
\\
I & a=L_{I}
\end{array}
\right.
\]
\edf


\bdf\label{df:fmlclasses}
\benu


\item
$A\in\Del_{0}$ {\rm iff there exists a} $\Del_{0}${\rm -formula} $\tht[\vec{x}]$ 
{\rm in the language} $\{\in\}$ {\rm and terms} $\vec{t}$ {\rm in} $\calL_{c}$
 {\rm such that} $A\equiv\tht[\vec{t}]$.
{\rm This means that} $A$ {\rm is bounded, and the predicates} $P,P_{I,n},Reg$ {\rm do not occur in} $A$.

\item
{\rm Putting} $\Sig_{0}:=\Pi_{0}:=\Del_{0}${\rm , the classes}
$\Sig_{m}$ {\rm and} $\Pi_{m}$ {\rm of formulae in the language} $\calL_{c}$
{\rm are defined as usual, where by definition}
$\Sig_{m}\cup\Pi_{m}\subset\Sig_{m+1}\cap\Pi_{m+1}$.

{\rm Each formula in} $\Sig_{m}\cup\Pi_{m}$ {\rm is in prenex normal form with alternating unbounded quantifiers and}
$\Del_{0}${\rm -matrix.}

 
 
 

\item
{\rm The set} $\Sig^{\Sig_{n+1}}(\lam)$ {\rm of sentences is defined recursively as follows. Let} $\{a,b,c\}\subset L_{I}$ {\rm and} $d\in L_{I}\cup\{L_{I}\}$.
 \benu
 \item
{\rm Each} $\Sig_{n+1}${\rm -sentence is in} $\Sig^{\Sig_{n+1}}(\lam)$.

 \item
 {\rm Each literal including} $Reg(a), P(a,b,c),P_{I,n}(a)$ {\rm and its negation
 is in} $\Sig^{\Sig_{n+1}}(\lam)$.

 \item
 $\Sig^{\Sig_{n+1}}(\lam)$ {\rm is closed under propositional connectives} $\lor,\land$.

 \item
{\rm Suppose}
 $\fal x\in d\, A(x)\not\in\Del_{0}${\rm . Then}
 $\fal x\in d\, A(x)\in \Sig^{\Sig_{n+1}}(\lam)$ {\rm iff} $A(\emptyset)\in \Sig^{\Sig_{n+1}}(\lam)$ {\rm and}
 $\rk_{L}(d)<\lam$.
 
 \item
{\rm Suppose}
 $\exi x\in d\, A(x)\not\in\Del_{0}${\rm . Then}
 $\exi x\in d\, A(x)\in \Sig^{\Sig_{n+1}}(\lam)$ {\rm iff} $A(\emptyset)\in \Sig^{\Sig_{n+1}}(\lam)$ {\rm and}
 $\rk_{L}(d)\leq\lam$.
 \eenu

\eenu
\edf

Note that the predicates $P,P_{I,n},Reg$ do not occur in $\Sig_{m}$-formulae.

\bdf\label{df:domFfml}
{\rm Let us extend the domain} $dom(F^{\Sig_{1}}_{x\cup\{\kap\}})=\mbox{{\rm Hull}}_{\Sig_{1}}^{I}(x\cup\{\kap\})$
{\rm of the Mostowski collapse to formulae.}
\[
dom(F^{\Sig_{1}}_{x\cup\{\kap\}})=\{A\in\Sig_{1}\cup\Pi_{1}: \sfk(A)\subset\mbox{{\rm Hull}}_{\Sig_{1}}^{I}(x\cup\{\kap\})\}
.\]
{\rm For} $A\in dom(F^{\Sig_{1}}_{x\cup\{\kap\}})$,
$F^{\Sig_{1}}_{x\cup\{\kap\}}" A$ {\rm denotes the result of replacing each constant} 
$c\in L_{I}$ {\rm by} 
$F^{\Sig_{1}}_{x\cup\{\kap\}}(c)${\rm , 
 each unbounded existential quantifier} $\exi z\in L_{I}$ {\rm by} 
 $\exi z\in L_{F^{\Sig_{1}}_{x\cup\{\kap\}}(I)}${\rm ,
and each unbounded universal quantifier} $\fal z\in L_{I}$ {\rm by} 
$\fal z\in L_{F^{\Sig_{1}}_{x\cup\{\kap\}}(I)}$.

{\rm For sequent, i.e., finite set of sentences} $\Gam\subset dom(F^{\Sig_{1}}_{x\cup\{\kap\}})${\rm , put}
 $F^{\Sig_{1}}_{x\cup\{\kap\}}"\Gam=\{F^{\Sig_{1}}_{x\cup\{\kap\}}" A: A\in\Gam\}$.
 
 {\rm Likewise the domain} $dom(F^{\Sig_{n}}_{x})=\mbox{{\rm Hull}}_{\Sig_{n}}^{I}(x)$
{\rm is extended to}
\[
dom(F^{\Sig_{n}}_{x})=\{A\in\Sig_{n}\cup\Pi_{n}: \sfk(A)\subset\mbox{{\rm Hull}}_{\Sig_{n}}^{I}(x)\}
\]
{\rm and for formula} $A\in dom(F^{\Sig_{n}}_{x})$,
$F^{\Sig_{n}}_{x}" A$, {\rm and sequent} $\Gam\subset dom(F^{\Sig_{n}}_{x})$,
 $F^{\Sig_{n}}_{x}"\Gam$
 {\rm are defined similarly.}
 
\edf

\bprp\label{prp:domFfml}
For $F=F^{\Sig_{1}}_{x\cup\{\kap\}}, F^{\Sig_{n}}_{x}$ and $A\in dom(F)$,
$A\lrarw F" A$.
\eprp

The assignment of disjunctions $A\simeq\bigvee(A_{\iota})_{\iota\in J}$ or
conjunctions  $A\simeq\bigwedge(A_{\iota})_{\iota\in J}$
to sentences $A$ is defined as in
\cite{Buchholz} {\it except\/} for $\Sig_{n}\cup\Pi_{n}$-sentences.

\bdf\label{df:assigndc}
\benu
\item\label{df:assigndc0}
{\rm If} $M$ {\rm is one of the literals} $a\in b,a\not\in b${\rm , then for} $J:=0$
\[
M:\simeq
\left\{
\begin{array}{ll}
\bigvee(A_{\iota})_{\iota\in J} & \mbox{{\rm if }} M \mbox{ {\rm is false (in }} $L$\mbox{{\rm )}}
\\
\bigwedge(A_{\iota})_{\iota\in J} &  \mbox{{\rm if }} M \mbox{ {\rm is true}}
\end{array}
\right.
\]

\item
$(A_{0}\lor A_{1}):\simeq\bigvee(A_{\iota})_{\iota\in J}$
{\rm and}
$(A_{0}\land A_{1}):\simeq\bigwedge(A_{\iota})_{\iota\in J}$
{\rm for} $J:=2$.

\item

$
Reg(a) :\simeq
\bigvee(a\not\in a)_{\iota\in J}
$
 {\rm and}
$
\lnot Reg(a):\simeq
\bigwedge(a\in a)_{\iota\in J} 
$
{\rm with}
$
J:=
\left\{
\begin{array}{ll}
1& \mbox{{\rm if }} a\in R
\\
0 & \mbox{{\rm otherwise}}
\end{array}
\right.
$.

\item
$
P(a,b,c) :\simeq
\bigvee(a\not\in a)_{\iota\in J} 
\mbox{ {\rm and }}
\lnot P(a,b,c) :\simeq
\bigwedge(a\in a)_{\iota\in J} 
$
{\rm with}
\[
J:=
\left\{
\begin{array}{ll}
1 & \mbox{{\rm if }}  
a\in R \spand \exi \alp<\ome_{n+1}(I+1)[b=\Psi_{a,n}\alp\spand
\alp\in\calh_{\alp,n}(b)\spand  c=F^{\Sig_{1}}_{b\cup\{a\}}(I)]
\\
0 & \mbox{{\rm otherwise}}
\end{array}
\right.
.\]

\item
$
P_{I,n}(a) :\simeq
\bigvee(a\not\in a)_{\iota\in J} 
\mbox{ {\rm and }}
\lnot P_{I,n}(a) :\simeq
\bigwedge(a\in a)_{\iota\in J} 
$
{\rm with}
\[
J:=
\left\{
\begin{array}{ll}
1 & \mbox{{\rm if }}  
 \exi \alp<\ome_{n+1}(I+1)[a=\Psi_{I,n}\alp\spand
\alp\in\calh_{\alp,n}(a)]
\\
0 & \mbox{{\rm otherwise}}
\end{array}
\right.
.\]

\item
{\rm Let} $\exi z\in b \, \tht[z]\in\Sig_{n}$ {\rm for} $b\in L_{I}\cup\{L_{I}\}$.
{\rm Then for the set}
\beqn\label{eq:dfmu}
d:=\mu z\in b \, \tht[z] :=
\min_{<_{L}}\{d : (d\in b \land \tht[d]) \lor (\lnot\exi z\in b\, \tht[z]\land d=0)\}
\eeqn
{\rm with a canonical well ordering} $<_{L}$ {\rm on} $L$, {\rm and}
$J=\{d\}$
\beqnarr
\exi z\in b\, \tht[z] & :\simeq & \bigvee(d\in b\land \tht[d])_{d\in J}
\label{eq:sigpimu}
\\
\fal z\in b\, \lnot\tht[z] & :\simeq & \bigwedge(d\in b \to \lnot\tht[d])_{d\in J}
\nonumber
\eeqnarr
{\rm where} $d\in b$ {\rm denotes a true literal, e.g.,} $d\not\in d$ {\rm when} $b=L_{I}$.


\item
{\rm Otherwise set for} $a\in L_{I}\cup\{L_{I}\}$ {\rm and} $J:=\{b: b\in a\}$
\[
\exi x\in a\, A(x):\simeq\bigvee(A(b))_{b\in J}
\mbox{ {\rm and }}
\fal x\in a\, A(x):\simeq\bigwedge(A(b))_{b\in J}
.\]

\eenu
\edf

Observe the following facts which are immediately seen from the definition.
Let $A\simeq\bigvee(A_{\iota})_{\iota\in J}$ or
 $A\simeq\bigwedge(A_{\iota})_{\iota\in J}$.
Then $\iota\in J$ is $Bool(\Sig_{n+1})$, and
 for $\iota\in J$, $x=A_{\iota}$ is a $\Del_{1}$-predicate $\vphi(A,\iota,x)$.

The rank $\rk(A)$ of sentences $A$ is defined by recursion on the  number of symbols occurring in 
$A$.

\bdf\label{df:rank}
\benu
\item\label{df:rank1}
$\rk(\lnot A):=\rk(A)$.

\item\label{df:rank2}
$\rk(a\in b):=0$.

\item\label{df:rank5-1}
$\rk(Reg(a)):=\rk(P(a,b,c)):=\rk(P_{I,n}(a)):=1$.


\item\label{df:rank6}
$\rk(A_{0}\lor A_{1}):=\max\{\rk(A_{0}),\rk(A_{1})\}+1$.


\item\label{df:rank7}
$
\rk(\exi x\in a\, A(x)):=
\max\{\ome\alp, \rk(A(\emptyset))+1\}
$
{\rm for} $\alp=\rk_{L}(a)$.
\eenu
\edf

\bprp\label{lem:rank}
Let $A\simeq\bigvee(A_{\iota})_{\iota\in J}$ or $A\simeq\bigwedge(A_{\iota})_{\iota\in J}$.
\benu
\item\label{lem:rank15-1}
$\fal\iota\in J(\sfk(A_{\iota})\subset\sfk(A)\cup\{\iota\})$.

\item\label{lem:rank15}
$A\in \Sig^{\Sig_{n+1}}(\lam)\Rarw\fal\iota\in J(A_{\iota}\in \Sig^{\Sig_{n+1}}(\lam))$.

\item\label{prp:rksig3}
For an ordinal $\lam\leq I$ with $\ome\lam=\lam$,
$
\rk(A)<\lam \Rarw A\in\Sig^{\Sig_{n+1}}(\lam)
$.


\item\label{lem:rank0}
$\rk(A)<I+\ome$.

\item\label{lem:rank1}
$\rk(A)\in\{\ome \, \rk_{L}(a)+i : 
a\in\sfqk(A)\cup\{I\}, 
i\in\ome\}
\subset\mbox{{\rm Hull}}_{\Sig_{1}}^{I}(\sfk(A))$.

\item\label{lem:rank2}
$\fal\iota\in J(\rk(A_{\iota})<\rk(A))$.

\eenu
\eprp
\bprf
\ref{lem:rank}.\ref{lem:rank2}.
This is seen from the fact that
$a\in b\in L_{I}\cup\{L_{I}\} \Rarw \rk_{L}(a)<\rk_{L}(b)$.
\eprf

\subsection{Operator controlled derivations}\label{subsec:opcontderivation}

Let $\calh$ be an operator $\calh:\calP(L_{I})\to\calP(L_{I}\cup\ome_{n+1}(I+1))$.
The map $X\mapsto\calh_{\alp,n}(X)$ defined in Definition \ref{df:Cpsiregularsm}
is an example of such an operator.
For $\Tht\in \calP(L_{I})$, $\calh[\Tht]$ denotes the operator defined by
$\calh[\Tht](X):=\calh(\Tht\cup X)$ for $X\in\calP(L_{I})$.

Let $\calh$ be an operator,
$\kap \in R^{+}$, $\Gam$ a sequent, $a<^{\veps}\ome_{n+1}(I+1)$ and 
$b<^{\veps}I\oplus\ome$.
We define a relation $(\calh,\kap,n)\vdash^{a}_{b}\Gam$, which is read `there exists an infinitary derivation
of $\Gam$ which is $(\kap,n)${\it -controlled\/} by $\calh$, and
 whose height is at most $a$ and its cut rank is less than $b$'.

Sequents are finite sets of sentences, and inference rules are formulated in one-sided sequent calculus.
In Definition \ref{df:controlderreg} let us write $a<b$ for $a<^{\veps}b$.

\bdf\label{df:controlderreg}
{\rm Define a relation} $(\calh,\kap,n)\vdash^{a}_{b}\Gam$ {\rm as follows.}

$(\calh,\kap,n)\vdash^{a}_{b}\Gam$ {\rm holds if} 
\beqn
\label{eq:controlder}
\{a\}\cup\sfk(\Gam)\subset\calh:=\calh(\emptyset)
\eeqn
{\rm and one of the following
cases holds:}


\bdes
\item[$(\bigvee)$]
$A\simeq\bigvee\{A_{\iota}: \iota\in J\}$, $A\in\Gam$ {\rm and there exist} $\iota\in J$
{\rm and}
 $a(\iota)<a$ {\rm such that}
\beqn\label{eq:bigveebnd}
\rk_{L}(\iota)<\kap \Rarw \rk_{L}(\iota)< a
\eeqn
{\rm and}
$(\calh,\kap,n)\vdash^{a(\iota)}_{b}\Gam,A_{\iota}$.

\item[$(\bigwedge)$]
$A\simeq\bigwedge\{A_{\iota}: \iota\in J\}$, $A\in\Gam$ {\rm and for every}
$\iota\in J$ {\rm there exists an} $a(\iota)<a$ 
{\rm such that}
$(\calh[\{\iota\}],\kap,n)\vdash^{a(\iota)}_{b}\Gam,A_{\iota}$.

\item[$(cut)$]
{\rm There exist} $a_{0}<a$ {\rm and} 
$C$
{\rm such that} $\mbox{{\rm rk}}(C)<b$
{\rm and}
$(\calh,\kap,n)\vdash^{a_{0}}_{b}\Gam,\lnot C$
{\rm and}
$(\calh,\kap,n)\vdash^{a_{0}}_{b}C,\Gam$.

\item[$(\mbox{P}_{\lam})$]
$\lam\in R$ {\rm and there exists} $\alp<\lam$ {\rm such that} 
$(\exi x,y<\lam[\alp<x \land P(\lam,x,y)])\in\Gam$.

\item[$(\mbox{F}^{\Sig_{1}}_{x\cup\{\lam\}})$] 

$\lam\in \calh\cap R$, $x=\Psi_{\lam,n}\bet\in\calh$ {\rm for a} $\bet$
{\rm and there exist} $a_{0}<a$, 
$\Gam_{0}\subset\Sig_{1}$ {\rm and} $\Lam$ {\rm such that}
$\sfk(\Gam_{0})\subset\mbox{{\rm Hull}}_{\Sig_{1}}^{I}((\calh\cap x)\cup\{\lam\})$,
$\Gam=\Lam\cup (F^{\Sig_{1}}_{x\cup\{\lam\}}"\Gam_{0})$
{\rm and}
$
(\calh,\kap,n)\vdash^{a_{0}}_{b}\Lam,\Gam_{0}$,
{\rm where} $F^{\Sig_{1}}_{x\cup\{\lam\}}$ {\rm denotes the Mostowski collapse}
$F^{\Sig_{1}}_{x\cup\{\lam\}}: \mbox{{\rm Hull}}^{I}_{\Sig_{1}}(x\cup\{\lam\})\lrarw L_{F^{\Sig_{1}}_{x\cup\{\lam\}}(I)}$.

\item[$(\mbox{P}_{I,n})$]
{\rm There exists} $\alp<I$ {\rm such that} 
$(\exi x<I[\alp<x \land P_{I,n}(x)])\in\Gam$.

\item[$(\mbox{F}^{\Sig_{n}}_{x})$] 
$
x=\Psi_{I,n}\bet\in\calh
$ {\rm for a} $\bet$ {\rm and there exist} $a_{0}<a$, $\Gam_{0}\subset\Sig_{n}$
{\rm and} $\Lam$ {\rm such that}
$\sfk(\Gam_{0})\subset\mbox{{\rm Hull}}_{\Sig_{n}}^{I}(\calh\cap x)$,
$\Gam=\Lam\cup (F^{\Sig_{n}}_{x}"\Gam_{0})$ {\rm and}
$(\calh,\kap,n)\vdash^{a_{0}}_{b}\Lam,\Gam_{0}$,
{\rm where} $F^{\Sig_{n}}_{x}$ {\rm denotes the Mostowski collapse}
$F^{\Sig_{n}}_{x}: \mbox{{\rm Hull}}^{I}_{\Sig_{n}}(x)\lrarw L_{F^{\Sig_{n}}_{x}(I)}$.


\edes
\edf

Since $x\in\calh_{\gam,n}[\Tht]$, $x=\Psi_{\lam,n}\bet$ and $\iota\in J$ in 
$A\simeq\bigvee(A_{\iota})_{\iota\in J}$ or $A\simeq\bigwedge(A_{\iota})_{\iota\in J}$
are all first-order definable,
there exists a first-order formula (roughly estimated a $\Pi_{n+2}$-formula) $H_{n}$
such that the relation $(\calh_{\gam,n}[\Tht],\kap,n)\vdash^{a}_{b}\Gam$
is a fixed point of $H_{n}$ as in (\ref{eq:fixH}).
An inspection to Definition \ref{df:controlderreg}
shows that $H_{n}$ is strictly positive.

In what follows the relation should be understood as a fixed point of $H_{n}$,
and recall that we are working in 
the intuitionistic fixed point theory $\mbox{FiX}^{i}(\sfZFL)$
over \sfZF+(V=L) defined in subsection \ref{subsec:intfixZFL}.

\bprp\label{prp:kapcontrolled}
$(\calh,\kap,n)\vdash^{a}_{b}\Gam \spand \lam\leq\kap \Rarw (\calh,\lam,n)\vdash^{a}_{b}\Gam$.
\eprp

The inferences rules $(\bigvee)$, $(\bigwedge)$ and $(cut)$ 
are standard except
$\Sig_{n}\cup\Pi_{n}$-formulae are derived from specific minor formulae, (\ref{eq:sigpimu}).
$(\mbox{{\bf P}}_{\lam})$ is an axiom for deducing the axiom (\ref{eq:Z2}), 
$Reg(\lam) \to \fal z<\lam(\exi x,y<\lam[z<x\land P(\lam,x,y)])$,
and $(\mbox{{\bf F}}^{\Sig_{1}}_{x\cup\{\lam\}})$ for proving the axiom (\ref{eq:Z1}),
$P(\lam,x,y)\land z<x \to \vphi[\lam,z] \to \vphi^{y}[x,z]$ for $\Sig_{1}$ $\vphi$.
Likewise $(\mbox{{\bf P}}_{I})$ and $(\mbox{{\bf F}}^{\Sig_{n}}_{x})$ for the axioms  (\ref{eq:Z5}) and  (\ref{eq:Z4}).

Let us explain the purpose of the unusual(, though correct) inference rules $(\bigvee)$, $(\bigwedge)$
for deriving $\Sig_{n}\cup\Pi_{n}$-formulae.
For simplicity set $\lam=\ome_{1}$ and $F_{x}=F^{\Sig_{1}}_{x\cup\{\ome_{1}\}}$, and consider the language of ordinals.
Consider the standard inference rules for introducing existential quantifiers in which any correct witness
can be a witness:
\[
\infer{(\calh,\kap,n)\vdash\exi z<\bet\,\tht[\vec{\gam},z],\Gam}{(\calh,\kap,n)\vdash\tht[\vec{\gam},\alp],\Gam}
\]
where $\alp<\bet$.
Then its dual should be
\[
\infer{(\calh,\kap,n)\vdash\Gam,\fal z<\bet\,\lnot\tht[\vec{\gam},z]}
{\{(\calh[\{\alp\}],\kap,n)\vdash\Gam,\lnot\tht[\vec{\gam},\alp]\}_{\alp<\bet}}
\]
But then, we have to examine all possible witnesses $\alp<\bet$ in deriving the axiom
$\fal z<I\lnot\tht[z,\ome_{1},a],\exi z<F_{x}(I)\,\tht[z,F_{x}(\ome_{1}),a]$ for $a<x=F_{x}(\ome_{1})$:
Assume $a,x,y\in\calh$.
\[
\infer{(\calh,\kap,n)\vdash \fal z<I\lnot\tht[z,\ome_{1},a],\exi z<F_{x}(I)\,\tht[z,F_{x}(\ome_{1}),a]}
{\{(\calh[\{\alp\}],\kap,n)\vdash \lnot\tht[\alp,\ome_{1},a],\exi z<F_{x}(I)\,\tht[z,F_{x}(\ome_{1}),a]\}_{\alp<I}}
\]
For $\alp\in dom(F_{x})$ we can deduce it by $(\mbox{{\bf F}}_{x})$
\[
\infer{(\calh[\{\alp\}],\kap,n)\vdash \lnot\tht[\alp,\ome_{1},a],\exi z<F_{x}(I)\,\tht[z,F_{x}(\ome_{1}),a]}
{
\infer[(\mbox{{\bf F}}_{x})]{(\calh[\{\alp\}],\kap,n)\vdash \lnot\tht[\alp,\ome_{1},a], \tht[F_{x}(\alp),F_{x}(\ome_{1}),a]}
{(\calh[\{\alp\}],\kap,n)\vdash \lnot\tht[\alp,\ome_{1},a], \tht[\alp,\ome_{1},a]}
}
\]
But there are ordinals $\alp<I$ such that $\alp\not\in dom(F_{x})$ since 
$dom(F_{x})=\mbox{{\rm Hull}}_{\Sig_{1}}^{I}(x\cup\{\ome_{1}\})$ is countable, and $I>\ome_{1}$ is uncountable.

Moreover the same trouble occurs, when an inference rule for quantifiers followed by an $(\mbox{{\bf F}}_{x})$:
\[
\infer[(\mbox{{\bf F}}_{x})]{\Gam,\exi z<F_{x}(\bet)\,\tht[F_{x}(\vec{\gam}),z]}
{
\infer{\Gam,\exi z<\bet\,\tht[\vec{\gam},z]}{\Gam,\tht[\vec{\gam},\alp]}
}
\]
Even if $\alp<\bet$, it may be the case $\alp\not\in dom(F_{x})$. 
Then one can not replace a cut inference with its cut formula $\exi z<F_{x}(\bet)\,\tht[F_{x}(\vec{\gam}),z]$
by one of a cut formula of the form $\tht[F_{x}(\vec{\gam}), F_{x}(\alp^{\prime})]$.

Contrary to this, in the inference rule for $\del=\mu z<\bet\,\tht[\vec{\gam},z]$,
\[
\infer[(\bigvee)]{(\calh,\kap,n)\vdash\Gam,\exi z<\bet\,\tht[\vec{\gam},z]}{(\calh,\kap,n)\vdash\Gam,\tht[\vec{\gam},\del]}
\]
$\del$ is $\Sig_{1}$-definable from $\{\bet\}\cup\vec{\gam}$ if $\bet<I$.
Therefore if $\{\bet\}\cup\vec{\gam}\subset dom(F_{x})$, then so is $\del$.
\\

We will state some lemmata for the operator controlled derivations with sketches of their proofs
since
these can be shown as in \cite{Buchholz}.

In what follows by an operator we mean an $\calh_{\gam,n}[\Tht]$ for a finite set $\Tht$ of sets.

\blem\label{lem:tautology}{\rm (Tautology)}
If $\sfk(\Gam\cup\{A\})\subset\calh$
then
$(\calh,I,n)\vdash^{2\footnotesize{\rk}(A)}_{0}\Gam,\lnot A, A$.
\elem

\blem\label{lem:completeness}{\rm ($\Sig_{n}\cup\Pi_{n}$-completeness)}
\\
For any sentence $A\in\Sig_{n}\cup\Pi_{n}$,
$
(A  \, \mbox{{\rm is true}}) \Rarw (\calh,I,n)\vdash^{2\footnotesize{\rk}(A)}_{0}A
$.
\elem
\bprf
This is seen by induction on the number of symbols occurring in $\Sig_{n}\cup\Pi_{n}$-sentences $A$.
\eprf

\blem\label{lem:elimfalse}{\rm (Elimination of false $\Sig_{n}$-sentences)}
\\
For any sentence $A\in\Sig_{n}$,
$
(A \, \mbox{{\rm is false}}) \spand (\calh,I,n)\vdash^{a}_{c}\Gam,A \Rarw 
(\calh,I,n)\vdash^{a}_{c}\Gam
$.
\elem
\bprf
This is seen by induction on $a$ using Proposition \ref{prp:domFfml}.
\eprf

\blem\label{lem:embedfund}
Let $\vphi[x,z]\in\Sig_{m}$ for $m\geq 1$, and $\Tht_{c}=\{\lnot\fal y(\fal x\in y\,\vphi[x,c]\to\vphi[y,c])\}$.
Then for any operator $\calh$, and any $a,c$,
$(\calh[\{c,a\}],I,n)\vdash^{I+2m+3+2\footnotesize{\rk}_{L}(a)}_{I+m+1}\Tht_{c},\fal x\in a\,\vphi[x,c]$.
\elem
{\bf Proof} by induction on $\rk_{L}(a)$.
Let $f(a)=I+2m+3+2\rk_{L}(a)$.
By IH we have for any $b\in a$,
$(\calh[\{c,b\}],I,n)\vdash^{f(b)}_{I+m+1}\Tht_{c},\fal x\in b\,\vphi[x,c]$.
On the other hand by Lemma \ref{lem:tautology} with 
$\rk(\vphi)\leq I+m-1$
 and $\rk(\fal x\in b\,\vphi[x,c])\leq I+m$,
we have $(\calh[\{c,b\}],I,n)\vdash^{I+2m+3}_{0}\Tht_{c},\lnot\fal x\in b\,\vphi[x,c], \vphi[b,c]$.
By a (cut) with $I+2m+3\leq f(b)$ we obtain
$(\calh[\{c,b\}],I,n)\vdash^{f(b)+1}_{I+m+1}\Tht_{c},\vphi[b,c]$.
$(\bigwedge)$ yields 
$(\calh[\{c,a\},I,n)\vdash^{f(a)}_{I+m+1}\Tht_{c},\fal x\in a\,\vphi[x,c]$.
\eprf

\bdf
$(\calh,I,n)\vdash_{c}^{<\alp}\Gam:\Lrarw \exi\bet<\alp[(\calh,I,n)\vdash_{c}^{\bet}\Gam]$.
\edf

\blem\label{th:embedreg}
Let $A$ be an axiom in $\mbox{{\rm T}}(I,n)$ except Foundation axiom schema.
Then
 $(\calh,I,n)\vdash_{0}^{<I+\ome}A$ for any operator $\calh=\calh_{\gam,n}$.
\elem
\bprf
By Lemma \ref{lem:completeness}
there remains nothing to show for $\Pi_{2}$-axioms in $\mbox{{\rm KP}}\ome+(V=L)$.

We consider the axiom (\ref{eq:Z1}).
Let a $\Sig_{1}$-formula $\vphi[x,a]\equiv\exi z\in L_{I}\,\tht[z,x,a]$ be given,
and assume $\lam,\iota,\nu,a\in L_{I}$.
\bdes
\item[Case 1]: $\lam\in R\spand \iota=\Psi_{\lam,n}\alp$ with 
$\alp\in\calh_{\alp,n}(\iota) \spand \nu=F^{\Sig_{1}}_{\iota\cup\{\lam\}}(I)$ and 
$a\in L_{\iota}$.

We show 
$
(\calh[\{\lam,\iota,a\}],I,n)\vdash^{<I}_{0} \lnot\vphi[\lam,a],\vphi^{\nu}[\iota,a]
$.

Let
$c=\mu z\in L_{I} \tht[z,\lam,a]$.
Then 
$\rk(\tht[c,\lam,a])<I$ since $\tht$ is $\Del_{0}$,
and by Lemma \ref{lem:tautology} we obtain
$
(\calh[\{\lam,a\}],I,n)\vdash^{<I}_{0} \lnot\tht[c,\lam,a],\tht[c,\lam,a]
$, where $c\in \mbox{{\rm Hull}}_{\Sig_{1}}^{I}(\{\lam,a\})\subset\calh[\{\lam,a\}]$.

By $(\mbox{{\bf F}}^{\Sig_{1}}_{\iota\cup\{\lam\}})$ with 
$\iota=F^{\Sig_{1}}_{\iota\cup\{\lam\}}(\lam)$, $a=F^{\Sig_{1}}_{\iota\cup\{\lam\}}(a)$,
$
(\calh[\{\lam,\iota,a\}],I,n)\vdash^{<I}_{0} 
\lnot\tht[c,\lam,a],\tht[F^{\Sig_{1}}_{\iota\cup\{\lam\}}(c),\iota,a]
$, where 
$F^{\Sig_{1}}_{\iota\cup\{\lam\}}(c)\in \calh[\{\lam,\iota,a\}]$
by $c \in \calh[\{\lam,\iota,a\}]$.

By $F^{\Sig_{1}}_{\iota\cup\{\lam\}}(c)=\mu z\in L_{\nu} \tht[z,\iota,a]\in L_{\nu}$
for $\nu=F^{\Sig_{1}}_{\iota\cup\{\lam\}}(I)$, and $(\bigvee)$,
$
(\calh[\{\lam,\iota,a\}], I,n)\vdash^{<I}_{0} \lnot\tht[c,\lam,a],\vphi^{\nu}[\iota,a]
$, where $\nu\in\calh[\{\lam,\iota,a\}]$.
By $(\bigwedge)$ we conclude
$
(\calh[\{\lam,\iota,a\}],I,n)\vdash^{<I}_{0} \lnot\vphi[\lam,a],\vphi^{\nu}[\iota,a]
$.
\item[Case 2]: Otherwise.

Then $\lnot P(\lam,\iota,\nu)\simeq\bigwedge\emptyset$ or 
$\lnot(a\in L_{\iota})\simeq\bigwedge\emptyset$.
Hence 
$(\calh[\{\lam,\iota,\nu,a\}],I,n)\vdash^{0}_{0}\lnot P(\lam,\iota,\nu), \lnot(a\in L_{\iota})$.

\edes
So in any case,
$(\calh[\{\lam,\iota,\nu,a\}],I,n)\vdash^{<I}_{0}
\lnot P(\lam,\iota,\nu), \lnot(a\in L_{\iota}), \lnot\vphi[\lam,a],\vphi^{\nu}[\iota,a]$.

By $(\bigvee)$ and $(\bigwedge)$ we obtain
$
(\calh,I,n)\vdash^{I}_{0} 
\fal \lam,a,x,y\in L_{I} \{ P(\lam,x,y)\to a\in L_{x} \to \vphi[\lam,a] \to \vphi^{y}[x,a] \}
$.
Note that $P(\lam,x,y)\to a\in L_{x} \to \vphi[\lam,a] \to \vphi^{y}[x,a] $ is not a $\Sig_{n}$-formula since the predicate $P$ occurs in it.

Likewise the axiom (\ref{eq:Z4}) is derived by $(\mbox{{\bf F}}^{\Sig_{n}}_{x})$, and
$
(\calh,I,n)\vdash^{<I+\ome}_{0} 
\mbox{(\ref{eq:Z4})}
$.

Finally consider the axiom (\ref{eq:Z3}).
If $a$ is not an ordinal, then 
$(\calh[\{a\}],I,n)\vdash^{<I}_{0}a\not\in Ord$
for a $\Del_{0}$-formula $Ord$.
Hence
$(\calh[\{a\}],I,n)\vdash^{<I}_{0}a\in Ord\to \exi y[y>a\land Reg(y)]$.
Assume $a$ is an ordinal.
By Proposition \ref{prp:clshull}.\ref{prp:clshull.4} and $n\geq 2$
we have $a^{+}\in\calh[\{a\}]$, and
$(\calh[\{a\}],I,n)\vdash^{<I}_{0}a^{+}>a\land Reg(a^{+})$, and
$(\calh[\{a\}],I,n)\vdash^{<I}_{0}a\in Ord\to \exi y[y>a\land Reg(y)]$.
Therefore by $(\bigwedge)$ we obtain
$(\calh,I,n)\vdash^{I}_{0}\fal x\in Ord \exi y[y>x\land Reg(y)]$.
\eprf

\blem\label{th:embedregthm}{\rm (Embedding)}\\
If $\mbox{{\rm T}}(I,n)\vdash \Gam[\vec{x}]$, there are $m,k<\ome$ such that
for any $\vec{a}\subset L_{I}$,
 $(\calh[\vec{a}],I,n)\vdash_{I+m}^{I\cdot 2+k}\Gam[\vec{a}]$ for any operator 
 $\calh=\calh_{\gam,n}$.
\elem
\bprf
By Lemma \ref{lem:embedfund}
we have
$(\calh,I,n)\vdash^{I\cdot 2}_{I+m+1}\fal u,z(\fal y(\fal x\in y\,\vphi[x,z]\to\vphi[y,z])\to \vphi[u,z])$
for $\vphi[x,z]\in\Sig_{m}$.
By Lemmata \ref{lem:tautology} and \ref{th:embedreg} it suffices to consider inference rules of logical connectives.

Suppose $(\calh[\{a,b\}],I,n)\vdash_{I+m}^{I\cdot 2+k}\Gam[a,b],\tht[a]$
and $(\calh[\{a,b\}],I,n)\vdash_{I+m}^{I\cdot 2+k}\Gam[a,b],a\in b$
for any $a\in L_{I}$, where we suppress parameters for simplicity.
We show for $A\equiv\exi z\in b\,\tht[z]$
\beqn\label{eq:embedregthm}
\fal a\in L_{I}\{(\calh[\{a,b\}],I,n)\vdash_{I+m}^{<I\cdot 2+\ome}\Gam[a,b],A\}
\eeqn

If $\exi z\in b\,\tht[z]\not\in\Sig_{n}$, then there is nothing to prove.
Assume $\exi z\in b\,\tht[z]\in\Sig_{n}$.
If $\exi z\in b\,\tht[z]$ is true (in $L$), then by Lemma \ref{lem:completeness}
we have $(\calh[\{b\}],I,n)\vdash^{2\footnotesize{\rk}(A)}_{0}A$, and hence (\ref{eq:embedregthm}) follows.

Otherwise $(a\not\in b)\lor\lnot\tht[a]$ is true.
If $\tht[a]$ is false, by Lemma \ref{lem:elimfalse}
we have $(\calh[\{a,b\}],I,n)\vdash_{I+m}^{I\cdot 2+k}\Gam[a,b]$, and hence (\ref{eq:embedregthm}).
Otherwise $a\in b$ is false, by Lemma \ref{lem:elimfalse}
we have $(\calh[\{a,b\}],I,n)\vdash_{I+m}^{I\cdot 2+k}\Gam[a,b]$, and hence (\ref{eq:embedregthm}).

Next assume that $a$ does not occur in $\Gam$.
Then wee can choose $a$ as we wish.
If $\exi  z\in b\,\tht[z]$ is true, then let
$a=\mu z\in b\,\tht[z]\in\mbox{Hull}_{\Sig_{1}}^{I}(\sfk(A))$. 
$(\calh[\{b\}],I,n)\vdash_{I+m}^{<I\cdot 2+\ome}\Gam[b],\tht[a]$ 
yields (\ref{eq:embedregthm}) 
by $(\bigvee)$ and $\rk_{L}(a)<I$.
Otherwise let $a=0$.
$(0\not\in b)\lor\lnot\tht[0]$ is true.
The rest is the same as above.

The case for introducing a universal quantifier is similar to the existential case.
\eprf

\bcor\label{cor:embed}
Assume $\mbox{{\rm T}}(I,n)\vdash \tht$ for a sentence $\tht$.
Let $m_{0}$ be a number such that $\vphi\in\Sig_{m_{0}}$ 
if an instance
$\fal u,z(\fal y(\fal x\in y\,\vphi[x,z]\to\vphi[y,z])\to \vphi[u,z])$
of Foundation axiom schema occurs in the given $\mbox{{\rm T}}(I,n)$-proof of $\tht$.

Then for $m=\max\{m_{0}+6,n+5\}$,
 $(\calh,I,n)\vdash_{I+m}^{<I\cdot 2+\ome}\vphi$ for any operator 
 $\calh=\calh_{\gam,n}$.
\ecor
\bprf
This is seen from the proof of Lemma \ref{th:embedregthm}, and
$\rk(\fal u,z(\fal y(\fal x\in y\,\vphi[x,z]\to\vphi[y,z])\to \vphi[u,z]))\leq I+m_{0}+5$
if $\vphi\in\Sig_{m_{0}}$ and
$\rk(A)\leq I+n+4$ for the universal closure $A$ of instances of axioms (\ref{eq:Z0})-(\ref{eq:Z5}) in 
$\mbox{{\rm T}}(I,n)$.
Specifically for $A\equiv(\fal x,a(P_{I,n}(x) \to a\in L_{x} \to \vphi[a] \to \vphi^{x}[a]))$ of (\ref{eq:Z4})
with $\vphi\in\Sig_{n}$,
we have $\rk(A)\leq I+n+4$.
\eprf

\blem\label{lem:inversionreg}{\rm (Inversion)}
\\
Let  $d=\mu z\in b\,\tht[\vec{c},z] $ for $\tht\in \Pi_{n-1}$.
Then
\[
(\calh,\kap,n)\vdash^{a}_{b}\Gam, \exi z\in b \,\tht[\vec{c},z] \Rarw
(\calh,\kap,n)\vdash^{a}_{b}\Gam, d\in b\land \tht[\vec{c},d]
\]
and
\[
(\calh,\kap,n)\vdash^{a}_{b}\Gam, \fal z\in b\,\lnot\tht[\vec{c},z] \Rarw
(\calh,\kap,n)\vdash^{a}_{b}\Gam, d\in b \to \lnot\tht[\vec{c},d]
\]
\elem
\bprf
Consider the case when $\tht\in\Del_{0}$ and
$(\exi z\in b\,\tht[\vec{c},z])\equiv 
(\exi z\in F^{\Sig_{1}}_{\iota\cup\{\lam\}}(b_{0})\,\tht[F^{\Sig_{1}}_{\iota\cup\{\lam\}}(\vec{c}_{0}),z])$
is a main formula of an $(\mbox{{\bf F}}^{\Sig_{1}}_{\iota\cup\{\lam\}})$ for an $\iota=\Psi_{\lam,n}\alp$.

We have $\{b_{0}\}\cup\vec{c}_{0}\subset dom(F^{\Sig_{1}}_{\iota\cup\{\lam\}})$.
Then $d=\mu z\in F^{\Sig_{1}}_{\iota\cup\{\lam\}}(b_{0})\,\tht[F^{\Sig_{1}}_{\iota\cup\{\lam\}}(\vec{c}_{0}),z]=F^{\Sig_{1}}_{\iota\cup\{\lam\}}(d_{0})$
for
$d_{0}=\mu z\in b_{0}\,\tht[\vec{c}_{0},z]\in dom(F^{\Sig_{1}}_{\iota\cup\{\lam\}})$.
Thus $d_{0}\in b_{0}\land\tht[\vec{c}_{0},d_{0}]$ is a minor formula with its main $d\in b\land \tht[\vec{c},d]$ of the 
$(\mbox{{\bf F}}^{\Sig_{1}}_{\iota\cup\{\lam\}})$.
\eprf
\\

\blem\label{lem:reduction}{\rm (Reduction)}\\
Let  $C\simeq\bigvee(C_{\iota})_{\iota\in J}$.
\benu
\item\label{lem:reduction1}
Suppose
$C\not\in\{\exi x<\lam\exi y<\lam[\alp<x \land P(\lam,x,y)]: \alp<\lam\in R\}\cup\{\exi x<I[\alp<x \land P_{I,n}(x)]:\alp<I\}$.
Then
\[
(\calh,\kap,n)\vdash^{a}_{c}\Del,\lnot C \spand (\calh,\kap,n)\vdash^{b}_{c}C,\Gam \spand \rk(C)\leq c
\Rarw
(\calh,\kap,n)\vdash^{a+b}_{c}\Del,\Gam
\]
\item\label{lem:reduction2}
Assume $C\equiv(\exi x<\lam\exi y<\lam[\alp<x \land P(\lam,x,y)])$ for an $\alp<\lam\in R$ and $\bet\in\calh_{\bet,n}$.
Then
\[
(\calh_{\bet,n},\kap,n)\vdash^{a}_{b}\Gam,\lnot C 
\Rarw
(\calh_{\bet+1,n},\kap,n)\vdash^{a}_{b}\Gam
\]

\item\label{lem:reduction3}
Assume $C\equiv(\exi x<I[\alp<x \land P_{I,n}(x)])$ for an $\alp<I$ and $\bet\in\calh_{\bet,n}$.
Then
\[
(\calh_{\bet,n},\kap,n)\vdash^{a}_{b}\Gam,\lnot C 
\Rarw
(\calh_{\bet+1,n},\kap,n)\vdash^{a}_{b}\Gam
\]
\eenu
\elem
\bprf
\\
\ref{lem:reduction}.\ref{lem:reduction1} by induction on $b<^{\veps}\ome_{n+1}(I+1)$, 
cf. Lemma \ref{lem:vepsfix}.

Consider the case when both $C$ and $\lnot C$ are main formulae.
First consider the case when
$C\equiv (F^{\Sig_{1}}_{\iota\cup\{\lam\}}"\vphi)$ is a main formula of an $(\mbox{{\bf F}}^{\Sig_{1}}_{\iota\cup\{\lam\}})$ with a $\vphi\in\Sig_{1}$, 
and $\lnot C$ is a main formula of a $(\bigwedge)$.
Let 
$\lnot C\equiv \lnot F^{\Sig_{1}}_{\iota\cup\{\lam\}}"\vphi\equiv
\fal z\in F^{\Sig_{1}}_{\iota\cup\{\lam\}}(e)\lnot\tht[F^{\Sig_{1}}_{\iota\cup\{\lam\}}(\vec{e}),z]$ with 
$\vec{e}\subset\mbox{{\rm Hull}}_{\Sig_{1}}^{I}((\calh\cap \iota)\cup\{\lam\})$,
$e\in\mbox{{\rm Hull}}_{\Sig_{1}}^{I}((\calh\cap \iota)\cup\{\lam\})\cup\{I\}$
and for the set
$
d=\mu z\in F^{\Sig_{1}}_{\iota\cup\{\lam\}}(e)\, \tht[F^{\Sig_{1}}_{\iota\cup\{\lam\}}(\vec{e}),z]
$
its minor formula is $\lnot\tht[F^{\Sig_{1}}_{\iota\cup\{\lam\}}(\vec{e}),d]$.


For any $z\in\mbox{{\rm Hull}}_{\Sig_{1}}^{I}(\iota\cup\{\lam\})$ we have
$
F_{\iota,\lam}"\tht[\vec{e},z] \Lrarw \tht[\vec{e},z]
$.

Now consider the set
$
d_{0}=\mu z\in e\, \tht[\vec{e},z]
$.
Then $d_{0}\in \mbox{{\rm Hull}}_{\Sig_{1}}^{I}(\iota\cup\{\lam\})=dom(F^{\Sig_{1}}_{\iota\cup\{\lam\}})$, and
$F^{\Sig_{1}}_{\iota\cup\{\lam\}}(d_{0})=d$.
Moreover by $\{e\}\cup\vec{e}\subset\calh$ we have $d_{0}\in\calh$.

By Lemma \ref{lem:inversionreg},
 inversion on the main formula $\fal z\in F^{\Sig_{1}}_{\iota\cup\{\lam\}}(e)\lnot \tht[F^{\Sig_{1}}_{\iota\cup\{\lam\}}(\vec{e}),z]$ of the $(\bigwedge)$, 
  we
get $(\calh,\kap,n)\vdash^{a}_{c}\Del,\lnot\tht[F^{\Sig_{1}}_{\iota\cup\{\lam\}}(\vec{e}),F^{\Sig_{1}}_{\iota\cup\{\lam\}}(d_{0})]$, 
and inversion on the minor formula $\exi z\in e\,\tht[\vec{e},z]$ of $(\mbox{{\bf F}}^{\Sig_{1}}_{\iota\cup\{\lam\}})$
we get $(\calh,\kap,n)\vdash^{b_{0}}_{c}\tht[\vec{e},d_{0}]$ for the $d_{0}\in e$,
and then by $(\mbox{{\bf F}}^{\Sig_{1}}_{\iota\cup\{\lam\}})$ go back to
$\lnot\tht[F^{\Sig_{1}}_{\iota\cup\{\lam\}}(\vec{e}),F^{\Sig_{1}}_{\iota\cup\{\lam\}}(d_{0})]$.

Transfer
{\small
\[
\infer[(cut)]{\Del,F^{\Sig_{1}}_{\iota\cup\{\lam\}}"\Gam,\Lam}
{
\Del,\fal z\in F^{\Sig_{1}}_{\iota\cup\{\lam\}}(e)\lnot \tht[F^{\Sig_{1}}_{\iota\cup\{\lam\}}(\vec{e}),z]
&
 \infer[(\mbox{{\bf F}}^{\Sig_{1}}_{\iota\cup\{\lam\}})]{\exi z\in F^{\Sig_{1}}_{\iota\cup\{\lam\}}(e)\tht[F^{\Sig_{1}}_{\iota\cup\{\lam\}}(\vec{e}),z],F^{\Sig_{1}}_{\iota\cup\{\lam\}}"\Gam,\Lam}
 {
   \exi z\in e\, \tht[\vec{e},z],\Gam,\Lam
  }
}
\]
}
to
{\small
\[
\infer[(cut)]{\Del,F^{\Sig_{1}}_{\iota\cup\{\lam\}}"\Gam,\Lam}
{
 \Del,\lnot \tht[F^{\Sig_{1}}_{\iota\cup\{\lam\}}(\vec{e}),F^{\Sig_{1}}_{\iota\cup\{\lam\}}(d_{0})]
&
  \infer[(\mbox{{\bf F}}^{\Sig_{1}}_{\iota\cup\{\lam\}})]{\tht[F^{\Sig_{1}}_{\iota\cup\{\lam\}}(\vec{e}),F^{\Sig_{1}}_{\iota\cup\{\lam\}}(d_{0})],F^{\Sig_{1}}_{\iota\cup\{\lam\}}"\Gam,\Lam}
  {\tht[\vec{e},d_{0}],\Gam,\Lam}
}
\]
}
Next consider the case $(\mbox{{\bf F}}_{\iota})$ vs. $(\mbox{{\bf F}}_{\iota_{1}})$ with $\iota_{1}>\iota$,
where $F_{\iota}=F^{\Sig_{1}}_{\iota,\lam}$ for some $\lam\in R$ or $F_{\iota}=F^{\Sig_{n}}_{\iota}$ with $\lam=I$,
and similarly for $F_{\iota_{1}}$.

Let $F_{\iota}"\vphi$ be a main formula of $(\mbox{{\bf F}}_{\iota})$, and 
$\lnot F_{\iota}"\vphi\equiv\lnot F_{\iota_{1}}"\tht$
a main formula of $(\mbox{{\bf F}}_{\iota_{1}})$.

Then 
by $\iota_{1}>\iota$ and Proposition \ref{prp:DeldefF} we have $F_{\iota}(I)<\iota_{1}$, and hence
$F_{\iota_{1}}"F_{\iota}"\vphi\equiv F_{\iota}"\vphi\equiv F_{\iota_{1}}"\tht$, i.e.,
$\tht\equiv F_{\iota}"\vphi$.
\[
\infer[(cut)]{\Lam,F_{\iota}"\Gam,\Lam_{1},F_{\iota_{1}}"\Gam}
{
\infer[(\mbox{{\bf F}}_{\iota})]{\Lam,F_{\iota}"\Gam,F_{\iota}"\vphi}{\Lam,\Gam,\vphi}
&
\infer[(\mbox{{\bf F}}_{\iota_{1}})]{\lnot F_{\iota}"\vphi,\Lam_{1},F_{\iota_{1}}"\Gam_{1}}{\lnot F_{\iota}"\vphi,\Lam_{1},\Gam_{1}}
}
\]
\ref{lem:reduction}.\ref{lem:reduction2}.
Suppose
$C\equiv (\exi x<\lam\exi y<\lam[\alp<x \land P(\lam,x,y)])$.
We have $(\calh_{\bet,n},\kap,n)\vdash^{a}_{b}\Gam,\lnot\exi x<\lam\exi y<\lam[\alp<x \land P(\lam,x,y)]$ 
with $\alp<\lam$.

Let $\iota=\Psi_{\lam,n}\bet$ and $\nu=F^{\Sig_{1}}_{\iota\cup\{\lam\}}(I)$.
Since $\alp\in\calh_{\bet,n}\cap\lam$,
we have $\alp<\Psi_{\lam,n}\bet=\iota$.
Moreover by $\bet\in\calh_{\bet,n}$ we have $\iota,\nu\in\calh_{\bet+1,n}$.
By inversion 
$
(\calh_{\bet+1,n},\kap,n)\vdash^{a}_{b}\Gam,\lnot[\alp<\iota\land P(\lam,\iota,\nu)]
$
and once again by inversion with $\lnot P(\lam,\iota,\nu)\simeq(\lam\in\lam)$ 
we have
$
(\calh_{\bet+1,n},\kap,n)\vdash^{a}_{b}\Gam,\alp\not<\iota,\lam\in\lam
$.
By eliminating the false sentences $\alp\not<\iota,\lam\in\lam$ we have
$
(\calh_{\bet+1,n},\kap,n)\vdash^{a}_{b}\Gam
$.
\\
\ref{lem:reduction}.\ref{lem:reduction3}. This is seen as in Lemma \ref{lem:reduction}.\ref{lem:reduction2}
by introducing the ordinal $\iota=\Psi_{I,n}\bet$.
\eprf




In the following Lemma \ref{lem:predcereg}, note that 
$\rk(\exi x<\lam\exi y<\lam[\alp<x \land P(\lam,x,y)])=\lam+1$ for $\alp<\lam\in R$, and
$\rk(\exi x<I[\alp<x \land P_{I,n}(x)])=I$.

\blem\label{lem:predcereg}{\rm (Predicative Cut-elimination)}
\benu
\item\label{lem:predcereg2}
$(\calh,\kap,n)\vdash^{b}_{c+\ome^{a}}\Gam
\spand 
[c,c+\ome^{a}[\cap (\{\lam+1:\lam\in R\}\cup\{I\})=\emptyset
\spand a\in\calh
\Rarw (\calh,\kap,n)\vdash^{\vphi ab}_{c}\Gam$.

\item\label{lem:predcereg4}
For $\lam\in R$, if $\ome^{b}<\ome_{n+1}(I+1)$,
$(\calh_{\gam,n},\kap,n)\vdash^{b}_{\lam+2}\Gam \spand \gam\in\calh_{\gam,n}
 \Rarw 
(\calh_{\gam+b,n},\kap,n)\vdash^{\ome^{b}}_{\lam+1}\Gam$.

\item\label{lem:predcereg5}
If $\ome^{b}<\ome_{n+1}(I+1)$,
$(\calh_{\gam,n},\kap,n)\vdash^{b}_{I+1}\Gam \spand \gam\in\calh_{\gam,n}
 \Rarw 
(\calh_{\gam+b,n},\kap,n)\vdash^{\ome^{b}}_{I}\Gam$.

\item\label{lem:predcereg6}
$(\calh_{\gam,n},\kap,n)\vdash^{b}_{c+\ome^{a}}\Gam
\spand 
[c,c+\ome^{a}[\cap R^{+}=\emptyset
\spand a\in\calh_{\gam,n}
\Rarw (\calh_{\gam+\vphi ab,n},\kap,n)\vdash^{\vphi ab}_{c}\Gam$.
\eenu
\elem
\bprf
\ref{lem:predcereg}.\ref{lem:predcereg6}.
This follows from Lemmata \ref{lem:predcereg}.\ref{lem:predcereg2}, \ref{lem:predcereg}.\ref{lem:predcereg4} and \ref{lem:predcereg}.\ref{lem:predcereg5} using the facts $\vphi ab\geq b$, and $a>0\Rarw \vphi 0(\vphi ab)=\vphi ab$.
\eprf

\bdf
{\rm For a formula} $\exi x\in d\, A(x)$ {\rm and ordinals} $\lam=\rk_{L}(d)\in R^{+}, \alp$,
$(\exi x\in d\, A)^{(\exi\lam\restrict\alp)}$ {\rm denotes the result of restricting the} outermost existential quantifier 
$\exi x\in d$ {\rm to} $\exi x\in L_{\alp}$,
$(\exi x\in d\, A)^{(\exi\lam\restrict\alp)}\equiv
(\exi x\in L_{\alp}\, A)$.
\edf

In what follows $F_{x,\lam}$ denotes $F^{\Sig_{1}}_{x,\lam}$ when $\lam\in R$, and $F^{\Sig_{n}}_{x}$ when $\lam=I$.

\blem\label{lem:boundednessreg}{\rm (Boundedness)}
Let $\lam\in R^{+}$, $C\equiv(\exi x\in d\, A)$ and $C\not\in\{\exi x<\lam \exi y<\lam[\alp<x \land P(\lam,x,y)]: \alp<\lam\in R\}\cup\{\exi x<I[\alp<x \land P_{I,n}(x)]:\alp<I\}$.
Assume that $\rk(C)=\lam=\rk_{L}(d)$.
\benu
\item\label{lem:boundednessregexi}
$
(\calh,\lam,n)\vdash^{a}_{c}\Lam, C \spand a\leq b\in\calh\cap\lam
\Rarw (\calh,\lam,n)\vdash^{a}_{c}\Lam,C^{(\exi\lam\restrict b)}
$.

\item\label{lem:boundednessregfal}
$
(\calh,\kap,n)\vdash^{a}_{c}\Lam,\lnot C \spand b\in\calh\cap\lam 
\Rarw (\calh,\kap,n)\vdash^{a}_{c}\Lam,\lnot (C^{(\exi\lam\restrict b)})
$.
\eenu
\elem
{\bf Proof} by induction on $a<^{\veps}\ome_{n+1}(I+1)$, cf. Lemma \ref{lem:vepsfix}.

Note that if a main formula $F_{\iota,\sig}"\vphi$ of an $(\mbox{{\bf F}}_{\iota,\sig})$ 
is in $\Sig^{\Sig_{n+1}}(\lam)$, then
either $\sig\leq\lam$ and there occurs no bounded quantifier
$Q x<\lam$ in $F_{\iota,\sig}"\vphi$, 
or $\sig>\iota>\lam$ and $(F_{\iota,\sig}"\vphi)^{(\exi\lam\restrict b)}\equiv F_{\iota,\sig}"\vphi$.

Let $C\simeq \bigvee(C_{\iota})_{\iota\in J}$ for $C_{\iota}\equiv A(\iota)$, 
and $(\calh,\lam,n)\vdash^{a(\iota)}_{c}\Lam,C,C_{\iota}$ with an $a(\iota)<a$ for an $\iota\in J=d$.
Otherwise $C^{(\exi\lam\restrict b)}\equiv C$ by the definition.
Then $C^{(\exi\lam\restrict b)}\simeq\bigvee(C{\iota})_{\iota\in J^{\prime}}$
where $J^{\prime}=L_{b}$.
By the condition (\ref{eq:bigveebnd}) we have $\rk_{L}(\iota)<\lam \Rarw \rk_{L}(\iota)< a\leq b$, and hence 
$\iota\in L_{b}=J^{\prime}$.
By IH we have  Lemmata \ref{lem:boundednessreg}.\ref{lem:boundednessregexi} 
and \ref{lem:boundednessreg}.\ref{lem:boundednessregfal}.
\eprf

\blem\label{th:Collapsingthmreg1}{\rm (Collapsing)}\\
Let $\lam\in R^{+}$ and $\sig\in R^{+}\cup\{\ome_{\alp}: \mbox{{\rm limit }} \alp<I\}$.

Suppose $\{\gam,\lam,\sig\}\subset\calh_{\gam,n}[\Tht]$ with $\fal\rho\geq\lam[\Tht\subset\calh_{\gam,n}(\Psi_{\rho,n}\gam)]$, and 
$
\Gam\subset\Sig^{\Sig_{n+1}}(\lam)
$.
Let
$\mu=\sig+1$ if $\sig\in R^{+}$.
Otherwise $\mu=\sig$ if $\sig=\ome_{\alp}$ for a limit $\alp<I$.
Then for $\hat{a}=\gam+\ome^{\sig+a}$ and $\kap=\max\{\sig,\lam\}$,
if $\hat{a}<^{\veps}\ome_{n+1}(I+1)$,
\[
 (\calh_{\gam,n}[\Tht],\kap,n)\vdash^{a}_{\mu}\Gam \Rarw
  (\calh_{\hat{a}+1,n}[\Tht],\lam,n)\vdash^{\Psi_{\lam,n}\hat{a}}_{\Psi_{\lam,n}\hat{a}}\Gam.
\]
\elem
{\bf Proof} by main induction on $\mu$ with subsidiary induction on $a<^{\veps}\ome_{n+1}(I+1)$,
cf. Lemma \ref{lem:vepsfix}.

First note that $\Psi_{\lam,n}\hat{a}\in\calh_{\hat{a}+1,n}[\Tht]=\calh_{\hat{a}+1,n}(\Tht)$
since 
$\hat{a}=\gam+\ome^{\sig+a}\in\calh_{\gam,n}[\Tht]\subset\calh_{\hat{a}+1,n}[\Tht]$
 by the assumption,
$\{\gam,\lam,\sig,a\}\subset\calh_{\gam,n}[\Tht]$.
 
Assume $(\calh_{\gam,n}[\Tht][\Lam],\kap,n)\vdash^{a_{0}}_{\mu}\Gam_{0}$ with 
$\fal\rho\geq\lam[\Lam\subset\calh_{\gam,n}(\Psi_{\rho,n}\gam)]$.
Then by $\gam\leq\hat{a}$, we have for any $\rho\geq\lam$,
$\hat{a_{0}}\in \calh_{\gam,n}[\Tht][\Lam]\subset \calh_{\gam,n}(\Psi_{\rho,n}\gam)\subset\calh_{\hat{a},n}(\Psi_{\rho,n}\hat{a})$.
This yields
that
\beqn\label{eq:collapsethm}
a_{0}<a \Rarw \fal\rho\geq\lam(\Psi_{\rho,n}\widehat{a_{0}}<\Psi_{\rho,n}\hat{a})
\eeqn

Second observe that $\sfk(\Gam)\subset\calh_{\gam,n}[\Tht]\subset\calh_{\hat{a}+1,n}[\Tht]$ by $\gam\leq\hat{a}+1$.
 
 Third we have
 \beqn\label{eq:Collapsingthm}
\fal\rho\geq\lam[\sfk(\Gam)\subset\calh_{\gam,n}(\Psi_{\rho,n}\gam)]
\eeqn 
{\bf Case 1}.
First consider the case: $\Gam\ni A\simeq\bigwedge\{A_{\iota} :\iota\in J\}$
\[
\infer[(\bigwedge)]{(\calh_{\gam,n}[\Tht],\kap,n)\vdash^{a}_{\mu}\Gam}
{
 \{(\calh_{\gam,n}[\Tht\cup\{\iota\}],\kap,n)\vdash^{a(\iota)}_{\mu}\Gam, A_{\iota}:\iota\in J\}
 }\]
where $a(\iota)<a$ for any $\iota\in J$.
We claim that
\beqn\label{eq:sigcollapsfal}
\fal\iota\in J\fal\rho\geq\lam(\iota\in\calh_{\gam,n}(\Psi_{\rho,n}\gam))
\eeqn
Consider the case when $A\equiv\fal x\in b\, \lnot A^{\prime}$.
There are two cases to consider.
First consider the case when 
$J=\{d\}$ for the set
$d=\mu x\in b\, A^{\prime}$.
Then 
$\iota=d=(\mu x\in b\,  A^{\prime})\in\mbox{{\rm Hull}}_{\Sig_{n}}^{I}(\sfk(A))\subset \calh_{\gam,n}(\Psi_{\rho,n}\gam)$
by (\ref{eq:Collapsingthm}).

Otherwise 
$\rk_{L}(b)<\lam$, i.e., $b\in L_{\lam}$.
Let $\rho\geq\lam$.
We have 
$b\in\sfk(A)\subset\calh_{\gam,n}[\Tht]\subset \calh_{\gam,n}(\Psi_{\rho,n}\gam)$.
Hence $b\in\calh_{\gam,n}(\Psi_{\rho,n}\gam)\cap L_{\rho}$.
Since $\calh_{\gam,n}(\Psi_{\rho,n}\gam)\cap\rho\subset\Psi_{\rho,n}\gam$
and $\Psi_{\rho,n}\gam$ is a multiplicative number,
we have 
$\calh_{\gam,n}(L_{\Psi_{\rho,n}\gam})\cap L_{\rho}=
\calh_{\gam,n}(\Psi_{\rho,n}\gam)\cap L_{\rho}\subset L_{\Psi_{\rho,n}\gam}$.
Therefore $\iota\in b\in L_{\Psi_{\rho,n}\gam}\subset \calh_{\gam,n}(\Psi_{\rho,n}\gam)$
as desired.

Hence (\ref{eq:sigcollapsfal}) was shown.

SIH  yields
\[
\infer[(\bigwedge)]
{(\calh_{\hat{a}+1,n}[\Tht],\lam,n)\vdash^{\Psi_{\lam,n}\hat{a}}_{\Psi_{\lam,n}\hat{a}}\Gam}
{
\{(\calh_{\widehat{a(\iota)}+1,n}[\Tht\cup\{\iota\}],\lam,n)\vdash^{\Psi_{\lam,n}\widehat{a(\iota)}}_{\Psi_{\lam,n}\widehat{a(\iota)}}
\Gam, A_{\iota}:\iota\in J\}
}
\]
for $\widehat{a(\iota)}=\gam+\ome^{\sig+a(\iota)}$, since $\Psi_{\lam,n}\widehat{a(\iota)}<\Psi_{\lam,n}\hat{a}$ by 
(\ref{eq:collapsethm}).
\\
{\bf Case 2}.
Next consider the case for an $A\simeq\bigvee\{A_{\iota}:\iota\in J\}\in\Gam$ and an $\iota\in J$
with $a(\iota)<a$ and $\rk_{L}(\iota)<\kap \Rarw \rk_{L}(\iota)<a$
\[
\infer[(\bigvee)]{(\calh_{\gam,n}[\Tht],\kap,n)\vdash^{a}_{\mu}\Gam}
{(\calh_{\gam,n}[\Tht],\kap,n)\vdash^{a(\iota)}_{\mu}\Gam,A_{\iota}}
\]
Assume $\rk_{L}(\iota)<\lam$. We show $\rk_{L}(\iota)<\Psi_{\lam,n}\hat{a}$.
By $\Psi_{\lam,n}\gam\leq\Psi_{\lam,n}\hat{a}$,
 it suffices to show $\rk_{L}(\iota)<\Psi_{\lam,n}\gam$.
 
Consider the case when $A\equiv\exi x\in b\, A^{\prime}$.
There are two cases to consider.
First consider the case when 
$J=\{d\}$ for the set
$d=\mu x\in b\, A^{\prime}$.
Then 
$\iota=d=(\mu x\in b\,  A^{\prime})\in\mbox{{\rm Hull}}_{\Sig_{n}}^{I}(\sfk(A))$, and 
$\rk_{L}(\iota)\in\mbox{{\rm Hull}}_{\Sig_{n}}^{I}(\sfk(A))\subset \calh_{\gam,n}(\Psi_{\lam,n}\gam)$
by (\ref{eq:Collapsingthm}).
If $\rk_{L}(\iota)<\lam$, then $\rk_{L}(\iota)\in \calh_{\gam,n}(\Psi_{\lam,n}\gam)\cap\lam\subset\Psi_{\lam,n}\gam$.

Otherwise
we have $J=b\in\sfk(A)\subset\calh_{\gam,n}[\Tht]$, and
we can assume that $\iota\in\sfk(A_{\iota})\subset\calh_{\gam,n}[\Tht]$.
Otherwise set $\iota=0$.
We have $\rk_{L}(\iota)<\rk_{L}(b)\leq\lam$, and $\rk_{L}(\iota)\in \calh_{\gam,n}(\Psi_{\lam,n}\gam)\cap\lam\subset\Psi_{\lam,n}\gam$.

 SIH yields for $\widehat{a(\iota)}=\gam+\ome^{\sig+a(\iota)}$
\[
\infer[(\bigvee)]
{(\calh_{\hat{a}+1,n}[\Tht],\lam,n)\vdash^{\Psi_{\lam,n}\hat{a}}_{\Psi_{\lam,n}\hat{a}}}
{
(\calh_{\widehat{a(\iota)}+1,n}[\Tht],\lam,n)\vdash^{\Psi_{\lam,n}\widehat{a(\iota)}}_{\Psi_{\lam,n}\widehat{a(\iota)}}\Gam, A_{\iota}
}
\]
{\bf Case 3}.
Third consider the case for an $a_{0}<a$ and a $C$ with $\mbox{rk}(C)<\mu$.
\[
\infer[(cut)]{(\calh_{\gam,n}[\Tht],\kap,n)\vdash^{a}_{\mu}\Gam}
{
(\calh_{\gam,n}[\Tht],\kap,n)\vdash^{a_{0}}_{\mu}\Gam,\lnot C
&
(\calh_{\gam,n}[\Tht],\kap,n)\vdash^{a_{0}}_{\mu}C,\Gam
}
\]
{\bf Case 3.1}. $\rk(C)<\lam$.

We have by (\ref{eq:Collapsingthm}) $\sfk(C)\subset\calh_{\gam,n}(\Psi_{\lam,n}\gam)$.
Proposition \ref{lem:rank}.\ref{lem:rank1} yields 
$\mbox{{\rm rk}}(C)\in\calh_{\gam,n}(\Psi_{\lam,n}\gam)\cap\lam\subset\Psi_{\lam,n}\gam\leq\Psi_{\lam,n}\hat{a}$.
By Proposition \ref{lem:rank}.\ref{prp:rksig3}
we see that $\{\lnot C,C\}\subset\Sig^{\Sig_{n+1}}(\lam)$.

SIH yields for $\widehat{a_{0}}=\gam+\ome^{\sig+a_{0}}$
\[
\infer[(cut)]{(\calh_{\hat{a}+1,n}[\Tht],\lam,n)\vdash^{\Psi_{\lam,n}\hat{a}}_{\Psi_{\lam,n}\hat{a}}\Gam}
{
(\calh_{\widehat{a_{0}}+1,n}[\Tht],\lam,n)\vdash^{\Psi_{\lam,n}\widehat{a_{0}}}_{\Psi_{\lam,n}\widehat{a_{0}}}\Gam,\lnot C
&
(\calh_{\widehat{a_{0}}+1,n}[\Tht],\lam,n)\vdash^{\Psi_{\lam,n}\widehat{a_{0}}
}_{\Psi_{\lam,n}\widehat{a_{0}}}C,\Gam
}
\]
{\bf Case 3.2}.
$\lam\leq\rk(C)<\mu$ and $\rk(C)\not\in R^{+}$.

Let $\pi:=\min\{\pi\in R^{+}:\pi>\rk(C)\}$.
We have $\pi\in R$ and $\pi\in\calh_{\gam,n}[\Tht]$ by $\rk(C)\in\calh_{\gam,n}[\Tht]$
 and  Proposition \ref{prp:clshull}.\ref{prp:clshull.4}.
 
Then $\lam\leq\rk(C)<\pi<\mu$, and hence $\{\lnot C,C\}\subset\Sig^{\Sig_{n+1}}(\pi)$
by Proposition \ref{lem:rank}.\ref{prp:rksig3}.
SIH with $\max\{\pi,\sig\}=\sig=\kap$ yields for $\widehat{a_{0}}=\gam+\ome^{\sig+a_{0}}$ and $\bet=\Psi_{\pi,n}\widehat{a_{0}}$,
$
(\calh_{\widehat{a_{0}}+1,n}[\Tht],\pi,n)\vdash^{\bet}_{\bet}\Gam,\lnot C
$
and
$
(\calh_{\widehat{a_{0}}+1,n}[\Tht],\pi,n)\vdash^{\bet}_{\bet}C,\Gam
$.

Let $\mu^{\prime}=\ome_{\alp}+1<\bet$ for $\pi=\ome_{\alp+1}$.
Then $\bet=\mu^{\prime}+\ome^{\bet}$ and $[\mu^{\prime},\mu^{\prime}+\ome^{\bet}[\cap R^{+}=\emptyset$.
Moreover $\rk(C)<\bet$.
By a $(cut)$
\[
\infer[(cut)]{(\calh_{\widehat{a_{0}}+1,n}[\Tht],\pi,n)\vdash^{\bet+1}_{\mu^{\prime}+\ome^{\bet}}\Gam}
{
(\calh_{\widehat{a_{0}}+1,n}[\Tht],\pi,n)\vdash^{\bet}_{\mu^{\prime}+\ome^{\bet}}\Gam,\lnot C
&
(\calh_{\widehat{a_{0}}+1,n}[\Tht],\pi,n)\vdash^{\bet}_{\mu^{\prime}+\ome^{\bet}}C,\Gam
}
\]
Predicative Cut-elimination \ref{lem:predcereg} yields 
\[
(\calh_{\widehat{a_{0}}+\vphi\bet(\bet+1),n}[\Tht],\pi,n)\vdash^{\vphi\bet(\bet+1)}_{\mu^{\prime}}\Gam
\]

We have $\mu^{\prime}<\mu$. 
MIH with $\max\{\lam,\mu^{\prime}\}<\pi$  and Proposition \ref{prp:kapcontrolled} yields 
\[
(\calh_{\widehat{a_{1}}+1,n}[\Tht],\lam,n)\vdash^{\Psi_{\lam,n}\widehat{a_{1}}}_{\Psi_{\lam,n}\widehat{a_{1}}}\Gam
\]
for $\widehat{a_{1}}=\widehat{a_{0}}+\vphi\bet(\bet+1)+\ome^{\ome_{\alp}+\vphi\bet(\bet+1)}
=\gam+\ome^{\sig+a_{0}}+\ome^{\ome_{\alp}+\vphi\bet(\bet+1)}<\gam+\ome^{\sig+a}=\hat{a}$
by $a_{0}<a$, $\ome_{\alp}<\sig$ and $\bet<\sig$ with a strongly critical $\sig$.
Thus $\Psi_{\lam,n}\widehat{a_{1}}<\Psi_{\lam,n}\hat{a}$ and
$(\calh_{\hat{a}+1,n}[\Tht],\lam,n)\vdash^{\Psi_{\lam,n}\hat{a}}_{\Psi_{\lam,n}\hat{a}}\Gam$.
\\
{\bf Case 3.3}.
$\lam\leq\rk(C)<\mu$ and $\pi:=\rk(C)\in R^{+}$.

Then $C\in\Sig^{\Sig_{n+1}}(\pi)$ and $\pi\leq\sig$.
Also $\pi\in\calh_{\gam,n}[\Tht]$.
$C$ is either a sentence $\exi x<I[\alp<x \land P_{I,n}(x)]$ with $\pi=I$, 
or a sentence $\exi x\in d\, A(x)$ with
 $\sfqk(A)<\pi=\rk_{L}(d)\leq I$.

In the first case we have $\kap=\sig=I$, and
$(\calh_{\gam+1,n}[\Tht],I,n)\vdash^{a_{0}}_{I+1}\Gam$ by Reduction \ref{lem:reduction}.\ref{lem:reduction3},
 and IH yields the lemma.

Consider the second case.
From the right uppersequent, 
SIH with $\max\{\pi,\sig\}=\sig=\kap$ yields for $\widehat{a_{0}}=\gam+\ome^{\sig+a_{0}}$ and $\bet_{0}=\Psi_{\pi,n}\widehat{a_{0}}\in\calh_{\widehat{a_{0}}+1,n}[\Tht]$
\[
(\calh_{\widehat{a_{0}}+1,n}[\Tht],\pi,n)\vdash^{\bet_{0}}_{\bet_{0}} C,\Gam
\]
Then by Boundedness \ref{lem:boundednessreg}.\ref{lem:boundednessregexi} and $\bet_{0}\in\calh_{\widehat{a_{0}}+1,n}[\Tht]$,
we have
\[
(\calh_{\widehat{a_{0}}+1,n}[\Tht],\pi,n)\vdash^{\bet_{0}}_{\bet_{0}}C^{(\exi\pi\restrict \bet_{0})},\Gam
\]
On the other hand we have by Boundedness \ref{lem:boundednessreg}.\ref{lem:boundednessregfal}
from the left uppersequent
\[
(\calh_{\widehat{a_{0}}+1,n}[\Tht],\pi,n)\vdash^{a_{0}}_{\mu}\Gam,\lnot (C^{(\exi\pi\restrict \bet_{0})})
\]
Moreover we have
$\lnot (C^{(\exi\pi\restrict \bet_{0})})\in\Sig^{\Sig_{n+1}}(\pi)$.
SIH yields for 
$\widehat{a_{0}}<\widehat{a_{1}}=\widehat{a_{0}}+1+\ome^{\sig+a_{0}}=\gam+\ome^{\sig+a_{0}}+1+\ome^{\sig+a_{0}}<\gam+\ome^{\sig+a}=\hat{a}$
and $\bet_{1}=\Psi_{\pi,n}\widehat{a_{1}}$
\[
(\calh_{\widehat{a_{1}}+1,n}[\Tht],\pi,n)\vdash^{\bet_{1}}_{\bet_{1}}
\Gam,\lnot C^{(\exi\pi\restrict \bet_{0})}
\]
Now we have $\widehat{a_{i}}\in \calh_{\widehat{a_{i}},n}(\Psi_{\pi,n}\hat{a})$ and $\widehat{a_{i}}<\hat{a}$ for $i<2$, and hence
$\bet_{0}=\Psi_{\pi,n}\widehat{a_{0}}<\bet_{1}=\Psi_{\pi,n}\widehat{a_{1}}<\Psi_{\pi,n}\hat{a}$. 
Therefore $\rk(C^{(\exi\pi\restrict \bet_{0})})<\bet_{1}<\Psi_{\pi,n}\hat{a}$.

Consequently
\[
\infer[(cut)]{(\calh_{\widehat{a_{1}}+1,n}[\Tht],\pi,n)\vdash^{\bet_{1}+1}_{\bet_{1}}\Gam}
{
(\calh_{\widehat{a_{1}}+1,n}[\Tht],\pi,n)\vdash^{\bet_{1}}_{\bet_{1}}\Gam,
\lnot C^{(\exi\pi\restrict \bet_{0})}
&
(\calh_{\widehat{a_{0}}+1,n}[\Tht],\pi,n)\vdash^{\bet_{0}}_{\bet_{0}}C^{(\exi\pi\restrict \bet_{0})},\Gam
}
\]
Let $(\alp,\mu^{\prime},\bet_{2})=(\alp,\ome_{\alp}+1,\bet_{1})$ if $\pi=\ome_{\alp+1}$, and $(\alp,\mu^{\prime},\bet_{2})=(\bet_{1},\bet_{1},0)=(\bet_{1},\ome_{\bet_{1}},0)$ if $\pi=I$.
Then $\bet_{1}\leq\mu^{\prime}+\ome^{\bet_{2}}$ and $[\mu^{\prime},\mu^{\prime}+\ome^{\bet_{2}}[\cap R^{+}=\emptyset$.

Predicative Cut-elimination \ref{lem:predcereg} yields 
\[
(\calh_{\widehat{a_{1}}+\vphi\bet_{2}(\bet_{1}+1),n}[\Tht],\pi,n)\vdash^{\vphi\bet_{2}(\bet_{1}+1)}_{\mu^{\prime}}\Gam
\]
We have $\mu^{\prime}<\mu$. 
MIH with $\max\{\lam,\mu^{\prime}\}\leq\pi$ yields 
\[
(\calh_{\widehat{a_{2}}+1,n}[\Tht],\lam,n)\vdash^{\Psi_{\lam,n}\widehat{a_{2}}}_{\Psi_{\lam,n}\widehat{a_{2}}}\Gam
\]
for $\widehat{a_{2}}=\widehat{a_{1}}+\vphi\bet_{2}(\bet_{1}+1)+\ome^{\ome_{\alp}+\vphi\bet_{2}(\bet_{1}+1)}
=\gam+\ome^{\sig+a_{0}}+\ome^{\sig+a_{0}}+\ome^{\ome_{\alp}+\vphi\bet_{2}(\bet_{1}+1)}<\gam+\ome^{\sig+a}=\hat{a}$
by $a_{0}<a$, $\ome_{\alp}<\sig$ and $\bet_{1}<\sig$ with a strongly critical $\sig$.
Thus $\Psi_{\lam,n}\widehat{a_{2}}<\Psi_{\lam,n}\hat{a}$ and
$(\calh_{\hat{a}+1,n},\lam,n)\vdash^{\Psi_{\lam,n}\hat{a}}_{\Psi_{\lam,n}\hat{a}}\Gam$.
\\
{\bf Case 4}. Fourth consider the case for an $a_{0}<a$
\[
\infer[(\mbox{{\bf F}})]{(\calh_{\gam,n}[\Tht],\kap,n)\vdash^{a}_{\mu}\Gam}
{
(\calh_{\gam,n}[\Tht],\kap,n)\vdash^{a_{0}}_{\mu}\Lam,\Gam_{0}
}
\]
where $\Gam=\Lam\cup F"\Gam_{0}$ and either  
$F=F^{\Sig_{1}}_{x\cup\{\rho\}}$, $\Gam_{0}\subset\Sig_{1}$ for some $x$ and $\rho$,
or 
$F=F^{\Sig_{n}}_{x}$, $\Gam_{0}\subset\Sig_{n}$ for an $x$.
Then $\Lam\cup\Gam_{0}\subset\Sig^{\Sig_{n+1}}(\lam)$.
SIH yields the lemma.
\eprf

\bcor\label{cor:completeCE}
Suppose 
$\Gam\subset\Sig^{\Sig_{n+1}}(\ome_{1})$.
Assume $(\calh_{0,n},I,n)\vdash^{I\cdot 2+k}_{I+m}\Gam$ for some $m,k<\ome$
such that $b=\ome_{m}(I\cdot 3+k)<\ome_{n+1}(I+1)$.
Let $\bet=\Psi_{\ome_{1},n}(b)$ and $c=\vphi\bet\bet$.
Then
$(\calh_{b+1,n},\ome_{1},n)\vdash^{c}_{0}\Gam$.
\ecor
\bprf
Let $(\calh_{0,n},I,n)\vdash^{I\cdot 2+k}_{I+m}\Gam$.
By Predicative Cut-elimination \ref{lem:predcereg}.\ref{lem:predcereg5}
we have $(\calh_{0,n},I,n)\vdash^{\ome_{m-1}(I\cdot 2+k)}_{I+1}\Gam$.
Collapsing \ref{th:Collapsingthmreg1} yields
$(\calh_{b+1,n},\ome_{1},n)\vdash^{\bet}_{\bet}\Gam$.
By Predicative Cut-elimination \ref{lem:predcereg}.\ref{lem:predcereg2}
we obtain
$(\calh_{b+1,n},\ome_{1},n)\vdash^{c}_{0}\Gam$.
\eprf



\bprp\label{prp:truthpreserve}
For each sentence $A$ in the language $\{\in\}\cup L_{I}$ the following holds.
\benu
\item\label{prp:truthpreserve1}
$A\simeq\bigvee(A_{\iota})_{\iota\in J} \Rarw \fal\iota\in J(A_{\iota} \mbox{ is an }
 \{\in\}\cup L_{I}\mbox{-sentence})$,
 and similarly for the case $A\simeq\bigwedge(A_{\iota})_{\iota\in J}$.

\item\label{prp:truthpreserve2}
$A\simeq\bigvee(A_{\iota})_{\iota\in J} \Lrarw 
(L_{I}\models A \Lrarw \exi\iota\in J(L_{I}\models A_{\iota}))$.

\item\label{prp:truthpreserve3}
$A\simeq\bigwedge(A_{\iota})_{\iota\in J} \Lrarw 
(L_{I}\models A \Lrarw \fal\iota\in J(L_{I}\models A_{\iota}))$.

\item\label{prp:truthpreserve4}
$(\calh,\ome_{1},n)\vdash^{\alp}_{0}\Gam
\spand \alp<^{\veps}\ome_{n+1}(I+1)
  \Rarw L_{I}\models\bigvee\Gam$.

\eenu
\eprp
\bprf
Propositions \ref{prp:truthpreserve}.\ref{prp:truthpreserve1}-\ref{prp:truthpreserve}.\ref{prp:truthpreserve3} are straightforward.

Proposition \ref{prp:truthpreserve}.\ref{prp:truthpreserve4} is proved by 
induction on $\alp<^{\veps}\ome_{n+1}(I+1)$ using Propositions \ref{prp:truthpreserve}.\ref{prp:truthpreserve1}-\ref{prp:truthpreserve}.\ref{prp:truthpreserve3}
and the fact that 
$(\mbox{{\bf F}}^{\Sig_{1}}_{x\cup\{\lam\}})$ and $(\mbox{{\bf F}}^{\Sig_{n}}_{x})$
are truth-preserving, that is to say if the upper sequent of these inferences is true,
then so is the lower sequent, cf. Proposition \ref{prp:domFfml}.
\eprf
\\

Observe that everything in this section is formalizable in $\mbox{FiX}^{i}(\sfZFL)$,
i.e., we need the excluded middle only for $\{\in,=\}$-formulas.

\section{{\bf Proof} of Theorem \ref{th:mainthZ}}\label{sect:proof}
For a sentence $\exi x\in L_{\ome_{1}}\vphi$ in the language $\{\in,\ome_{1}\}$,
assume
$\sfZF+(V=L)\vdash\exi x\in L_{\ome_{1}}\,\vphi$.
Let $n_{0}\geq 2$ be the number such that in the given $\sfZF+(V=L)$-proof
instances of axiom schemata of Separation and Collection are
$\Sig_{n_{0}}$-Separation and $\Sig_{n_{0}}$-Collection,
and let $n_{1}$ the number such that in the given $\sfZF+(V=L)$-proof
instances of Foundation axiom schema are
applied to $\Sig_{n_{1}}$-formulae.
Let $m=\max\{n_{1}+6,n_{0}+5\}$, and let $n=m+1$.
Then by Lemma \ref{lem:regularset} and Corollary \ref{cor:embed}
we see that
the fact
$(\calh_{0,n},I,n)\vdash^{<I\cdot 2+\ome}_{I+m}\exi x\in L_{\ome_{1}}\vphi$ is provable in 
$\mbox{FiX}^{i}(\sfZFL)$.
We have $b=\ome_{m}(I\cdot 3+\ome)<^{\veps}\ome_{n}(I+1)$.
In what follows work in $\mbox{FiX}^{i}(\sfZFL)$.
Corollary \ref{cor:completeCE}
yields
$(\calh_{b+1,n},\ome_{1},n)\vdash^{c}_{0}\exi x\in L_{\ome_{1}}\vphi$
for
$\bet=\Psi_{\ome_{1},n}(b)$ and $c=\vphi\bet\bet$.
Boundedness \ref{lem:boundednessreg} yields
$(\calh_{b+1,n},\ome_{1},n)\vdash^{c}_{0}\exi x\in L_{c}\, \vphi$.
Then by Proposition \ref{prp:truthpreserve}.\ref{prp:truthpreserve4}
with $c<^{\veps}\Psi_{\ome_{1},n}\ome_{n}(I+1)$
we obtain
$\exi x\in L_{\Psi_{\ome_{1},n}\ome_{n}(I+1)} \vphi$.

Since the whole proof is formalizable in $\mbox{FiX}^{i}(\sfZFL)$, 
we conclude
$\mbox{FiX}^{i}(\sfZFL)\vdash  \exi x\in L_{\Psi_{\ome_{1},n}\ome_{n}(I+1)} \vphi$.
Finally Theorem \ref{th:consvintfix} yields
$\sfZF+(V=L)\vdash  \exi x\in L_{\Psi_{\ome_{1},n}\ome_{n}(I+1)} \vphi$.
This completes a proof of Theorem \ref{th:mainthZ}.
\\

\noindent
{\bf Remark}.
Using notation systems of infinitary derivations as in \cite{BuchholzMLQ}, 
it is reasonable to expect the following:

Over a weak base theory T,
$\sfZF+(V=L)$ is a conservative extension of
$\mbox{{\rm T}}+(V=L)+\{\exi x<\ome_{1}[x=\Psi_{\ome_{1},n}\ome_{n}(I+1)] : n<\ome\}$ with respect to a class of formulae depending on T.
\\



\noindent
Since any cut-free derivation of a first-order sentence is finite in depth, we actually have the following Corollary \ref{cor:finiteness}.
\bcor\label{cor:finiteness}
Assume $\sfZF+(V=L)\vdash \exi x<\ome\,\vphi$.
Then there exist $n,h<\ome$ such that
\[
(\calh_{\ome_{n}(I+1)+1,n},\ome_{1},n)\vdash^{h}_{0}
\exi x<\ome\, \vphi
.\]
\ecor
{\bf Problem}.
Let $g$ be the G\"odel number of a $\mbox{{\rm T}}(I)$-proof of $\exi x<\ome\,\vphi$,
and $h=H(g)$ a bound of depth of cut-free derivation. Note here that a number $n<\ome$ such that
$(\calh_{\ome_{n}(I+1)+1,n},\ome_{1},n)\vdash^{h}_{0}\exi x<\ome\, \vphi$ is calculable from $g$.
Then the map $H$ on $\ome$ seems not to be provably total in $\sfZF+(V=L)$,
i.e., $\sfZF+(V=L)\not\vdash\fal g\in\ome\exi h\in\ome[h=H(g)]$,
and $H\not\in L_{\Psi_{\ome_{1}}\veps_{I+1}}$.


The problem is to find a reasonable hierarchy of reals$\in{}^{\ome}\ome$ indexed by countable ordinals,
and to show that $H$ is too rapidly growing to be provably total in $\sfZF+(V=L)$.




\end{document}